\documentclass[11pt]{article}

\usepackage[margin=1in]{geometry}
\usepackage{amsmath,amssymb,amsthm,amsfonts}
\usepackage{mathtools}
\usepackage{enumitem}
\usepackage{graphicx}
\graphicspath{{figures/}}
\usepackage{booktabs}
\usepackage[hidelinks]{hyperref}

\theoremstyle{plain}
\newtheorem{theorem}{Theorem}[section]
\newtheorem{lemma}[theorem]{Lemma}
\newtheorem{proposition}[theorem]{Proposition}
\newtheorem{corollary}[theorem]{Corollary}

\theoremstyle{definition}
\newtheorem{definition}[theorem]{Definition}
\newtheorem{assumption}[theorem]{Assumption}
\newtheorem{condition}[theorem]{Condition}
\newtheorem{example}[theorem]{Example}

\theoremstyle{remark}
\newtheorem{remark}[theorem]{Remark}

\newcommand{\R}{\mathbb{R}}
\newcommand{\Rp}{\mathbb{R}_{+}}
\newcommand{\N}{\mathcal{N}}
\newcommand{\G}{\mathcal{G}}
\newcommand{\Ical}{\mathcal{I}}
\newcommand{\Ucal}{\mathcal{U}}
\newcommand{\Ecal}{\mathcal{E}}
\newcommand{\Fcal}{\mathcal{F}}
\newcommand{\Vcal}{\mathcal{V}}
\newcommand{\Om}{\Omega}
\newcommand{\one}{\mathbf{1}}
\newcommand{\ee}{\mathrm{e}}
\newcommand{\dd}{\mathrm{d}}
\newcommand{\tr}{\operatorname{tr}}
\newcommand{\Rz}{\mathcal{R}_{0}}
\newcommand{\doi}[1]{\href{https://doi.org/#1}{doi:#1}}
\DeclareMathOperator{\spec}{\sigma} \DeclareMathOperator{\dist}{dist}
\usepackage{titlesec}

\titlespacing*{\section}
  {0pt}%
  {1em}%
  {0.3em}%

\titlespacing*{\subsection}
  {0pt}
  {1.2em}
  {0.4em}

\titlespacing*{\subsubsection}
  {0pt}
  {1.0em}
  {0.3em}

\usepackage{authblk}

\title{\bfseries
	Structure-Preserving Universal Differential Equations for Epidemic Dynamics
}

\author[1]{Abdallah Alsammani\thanks{
		Corresponding author: \texttt{aalsammani@desu.edu}}}
\author[2]{Maia Martcheva}

\affil[1]{Department of Mathematical Sciences,
	Delaware State University, Dover, DE 19901, USA}

\affil[2]{Department of Mathematics,
	University of Florida, Gainesville, FL 32611, USA}

\date{}
\begin{document}
\maketitle

\begin{abstract}
	\noindent
	Universal differential equations (UDEs) combine mechanistic epidemic models with learned neural closures, but generic additive corrections may violate essential epidemiological properties. We introduce a class of structure-preserving epidemic UDEs in which learned closures parameterize intercompartmental fluxes and the force of infection factorizes through infected compartments. This architecture guarantees positivity, population boundedness, and invariance of the disease-free set. We derive a rank-one next-generation matrix and a closed-form basic reproduction number determined by the learned rates at the disease-free equilibrium. Under explicit conditions, we establish results on threshold dynamics, disease-free and endemic equilibria, stability, and uniform persistence. For a neural SEIR model, center manifold analysis characterizes the direction of bifurcation through derivatives of the learned transmission rate. We further establish fundamental identifiability limitations: closures that agree on observed epidemic trajectories can yield different reproduction numbers and bifurcation types. Quantitative bounds relate uncertainty in the reproduction number to closure regularity and the distance of observations from the disease-free equilibrium, while monotonicity provides one-sided bounds. A numerical example involving three closures illustrates these results. Our findings provide a rigorous framework for epidemiologically consistent scientific machine learning and clarify which dynamical properties can be recovered from epidemic data.
\end{abstract}

\noindent{\sloppy\textbf{Keywords:} Universal differential equations; scientific machine learning; epidemic dynamics; structure-preserving systems; basic reproduction number; dynamical systems; identifiability.\par}

\noindent{\sloppy\textbf{AMS subject classifications:} 92D30; 34D23; 37N25; 34C23; 68T07.\par}

\medskip

\section{Introduction}\label{sec:intro}

Compartmental epidemic models owe their interpretability to the
fact that every term encodes a hypothesized mechanism, and their
mathematical theory is correspondingly well developed
\cite{KermackMcKendrick1927,AndersonMay1991,Hethcote2000,
Thieme2003,Martcheva2015,BrauerCastilloChavezFeng2019}.
The same rigidity is a liability when a mechanism is unknown
or misspecified, as with behavioral feedback, heterogeneous
mixing, waning immunity, or time-varying interventions, since
the assumed functional form propagates into every downstream
inference. For example, adaptive risk awareness and behavioral
fatigue can influence epidemic trajectories and long-term
disease dynamics \cite{MohammedAlsammani2026}.

Universal differential equations \cite{Rackauckas2020UDE} relax this rigidity while retaining mechanism: the known part of the vector field is kept exactly and a neural network is inserted where the model is uncertain. Epidemiological applications appeared quickly \cite{Dandekar2020,Linka2022,Shaier2022}, but they are computational, and the properties that make a compartmental model epidemiologically meaningful are dynamical rather than statistical. A model that reproduces incidence data while permitting negative prevalence, unbounded population, or spontaneous outbreaks from an uninfected state is not an epidemic model. Moreover, the invariants practitioners extract from a fitted model, above all the basic reproduction number and the stability type of the disease-free state, presuppose an invariant disease-free manifold and a splitting of the infected subsystem into new infections and transitions \cite{DiekmannHeesterbeekMetz1990,vandenDriesscheWatmough2002}; an unconstrained additive closure destroys both.

This paper asks which UDE models retain the structure of an epidemic model, and what can be proved about them. Two principles suffice. First, the network should parameterize fluxes between compartments rather than derivatives of state variables, which makes positivity and population balance structural. Second, the learned force of infection should have the form $\lambda_\phi(x)=\sum_{j}\beta^\phi_j(x)x_j$, so that the disease-free set is invariant for every weight vector, before and during training, while the dependence of the transmission coefficients on the state remains unrestricted. Within the resulting admissible class we obtain a qualitative theory in which the classical hypotheses of nonlinear incidence theory \cite{CapassoSerio1978,LiuLevinIwasa1986,Korobeinikov2007} become architectural constraints that can be imposed or certified on a compact region.

The analytical tools are classical: quasi-positivity, next-generation matrices and M-matrix arguments \cite{DiekmannHeesterbeekMetz1990,vandenDriesscheWatmough2002}, Lyapunov functions and comparison, persistence theory, and center manifold reduction. What is new is what they yield for a learned closure. The reproduction number of an epidemic UDE depends on the learned closure only through its values at the disease-free equilibrium, which may lie outside the region of state space sampled by outbreak data; for a fixed transfer matrix it is a nonnegative linear functional of the transmission value there (Theorem \ref{thm:R0closed}). Theorem \ref{thm:nonid} turns this into a nonidentifiability statement of a specific kind: closures that coincide on the region of state space visited by an observed trajectory generate identical trajectories, for every observation operator, yet may have any reproduction number in the admissible range. Theorem \ref{thm:quantid} quantifies the phenomenon, bounding the error in the threshold by the Lipschitz constant of the closure times the distance of the data from the disease-free state, and Proposition \ref{prop:monoident} shows which monotone architectures yield one-sided identifiable bounds. Theorem \ref{thm:backward} expresses the direction of the bifurcation at threshold for a neural SEIR model through the first derivatives of the learned rate at the disease-free equilibrium, quantities computable from a trained network by automatic differentiation; Corollary \ref{cor:bifnonid} shows that these derivatives, and with them the bifurcation type, are as unconstrained by trajectory data as the threshold itself, even for closures with the same reproduction number. Section \ref{sec:numerics} computes an example of three closures that share one trajectory and differ in threshold and bifurcation type. The message is that agreement with trajectory data, however accurate, does not by itself establish the dynamical invariants on which epidemiological interpretation rests, and that the architecture and the data design near the disease-free state decide what is recoverable.

Section \ref{sec:back} reviews the background. Section \ref{sec:model} defines the admissible class and summarizes which assumptions each result requires. Section \ref{sec:wp} treats well-posedness, positivity, boundedness, and invariance; Section \ref{sec:R0} the disease-free equilibrium and the reproduction number; Section \ref{sec:endemic} endemic equilibria; Section \ref{sec:stability} local and global stability; Section \ref{sec:persistence} persistence and extinction; Section \ref{sec:bifurcation} bifurcation; Section \ref{sec:sensitivity} sensitivity; Section \ref{sec:identifiability} approximation and identifiability; Section \ref{sec:numerics} the numerical illustration. Sections \ref{sec:discussion} and \ref{sec:conclusion} discuss and conclude.

\section{Background}\label{sec:back}

The qualitative analysis of compartmental epidemic models is built on a small set of structural principles that determine their long-term dynamics. Positivity of solutions follows from the quasi-positivity of the vector field \cite{Walter1998,Smith1995}, while demographic balance yields a linear differential inequality that defines a compact absorbing region. Threshold dynamics are characterized by the spectral radius of the next-generation operator \cite{DiekmannHeesterbeekMetz1990,vandenDriesscheWatmough2002,DiekmannHeesterbeekRoberts2010}, whose relationship to local stability is established through M-matrix theory \cite{BermanPlemmons1994}. Global stability is typically obtained using Volterra-type Lyapunov functions together with the LaSalle invariance principle \cite{KorobeinikovMaini2004,Korobeinikov2007,LaSalle1976}, through the matrix-theoretic and graph-theoretic constructions of Shuai and van den Driessche \cite{ShuaiVdD2013}, or through the geometric approach of Li and Muldowney \cite{LiMuldowney1996}. Persistence theory characterizes the dynamics above the epidemic threshold \cite{Thieme1993,FreedmanRuanTang1994,HaleWaltman1989,SmithThieme2011,Zhao2017}, whereas center manifold analysis distinguishes forward and backward bifurcations near criticality \cite{Carr1981,CastilloChavezSong2004,HadelerVandenDriessche1997,Dushoff1998,Brauer2004}. Collectively, these results show that the qualitative behavior of epidemic models depends strongly on the structural properties of the incidence function, not only on parameter values. Saturating \cite{CapassoSerio1978}, power-law \cite{LiuLevinIwasa1986}, and nonmonotone incidence \cite{RuanWang2003,XiaoRuan2007} produce different dynamical regimes, whereas monotonicity and sublinearity guarantee global stability for broad classes of nonlinear incidence functions \cite{Korobeinikov2007}.

On the learning side, neural differential equations represent the vector field with a neural network \cite{Chen2018NeuralODE,Kidger2022}, linking deep learning with dynamical systems \cite{WeinanE2017,HaberRuthotto2017}; their representational limitations motivated augmented architectures and structural analyses \cite{Dupont2019,Massaroli2020}. Physics-informed neural networks incorporate the residual of a governing equation into the loss function \cite{Lagaris1998,RaissiPerdikarisKarniadakis2019,Karniadakis2021}, supported by a growing error theory \cite{MishraMolinaro2022,DeRyckMishra2022}; sparse regression identifies governing equations from candidate libraries \cite{BruntonProctorKutz2016}; and universal differential equations retain the mechanistic model while learning only its unknown components \cite{Rackauckas2020UDE,RackauckasNie2017}. Universal approximation in $C^1$ on compact sets \cite{Cybenko1989,Hornik1991,HornikStinchcombeWhite1990,Pinkus1999} ensures that such closures are sufficiently expressive for the framework developed here. Positivity of learned dynamics has recently been enforced at the level of the numerical scheme by parameterizing production and destruction rates within a nonstandard finite difference discretization \cite{Zinihi2026}, and the question of which unknown functional components of a hybrid model can be recovered from data has been formalized as functional identifiability \cite{LomanBaker2025,LomanBrowningBaker2026}.

Three gaps motivate this work. First, although positivity can be preserved by structured parameterizations of the learned rates \cite{Zinihi2026}, we are not aware of a qualitative theory of epidemic UDEs: nothing in the additive formulation guarantees that the disease-free set remains invariant, that a threshold exists, or that the stability and bifurcation structure of the mechanistic model survives the insertion of a network. In the standard additive formulation,
\[
\dot{x}=f(x,\theta)+\N_\phi(x),
\]
where $\N_\phi:\R^n\to\R^n$ is a neural-network correction to the mechanistic vector field, parameterized by weights $\phi$. Following the usage of turbulence and reduced-order modeling, such a learned term is called a \emph{closure}: it stands in for processes that the mechanistic model leaves unresolved, here an unknown transmission or transition mechanism, and in Section \ref{sec:model} it is restricted to act through intercompartmental fluxes, which distinguishes it from an arbitrary additive perturbation. In the additive form, generic neural perturbations need not preserve quasi-positivity or the disease-free set, and violations are typically addressed through clipping or penalty methods, which provide numerical rather than structural remedies. Second, we are not aware of a general framework for defining the basic reproduction number of an epidemic UDE; a threshold computed from the mechanistic component alone ignores the learned closure. Third, functional identifiability as developed in \cite{LomanBaker2025,LomanBrowningBaker2026} asks whether an unknown function is determined by the outputs generated from all initial conditions; what a single observed epidemic determines about dynamical invariants such as the threshold and the bifurcation type does not appear to have been addressed, although this is an intrinsic property of the inverse problem rather than of the optimization algorithm.

\section{The admissible class}\label{sec:model}

\subsection{Mechanistic model and standing assumptions}

Let $n\ge2$ and let $x=(x_1,\dots,x_n)^{\!\top}\in\Rp^n$ be the compartment vector of a population with $n$ epidemiological classes, evolving under
\begin{equation}\label{eq:mech}
\dot x=f(x,\theta),\qquad x(0)=x_0\in\Rp^n,
\end{equation}
with $\theta$ in a compact parameter set $\Theta\subset\R^p$. Fix a partition $\{1,\dots,n\}=\Ical\cup\Ucal$ into infected classes $\Ical$, $|\Ical|=m\ge1$, and uninfected classes $\Ucal$, $|\Ucal|=n-m\ge1$, and denote by $y=(x_i)_{i\in\Ical}$ and $z=(x_i)_{i\in\Ucal}$ the infected and uninfected coordinates. We identify $x$ with the pair $(z,y)$, so that $(z,0)$ denotes the point with uninfected coordinates $z$ and no infection; in block matrices, such as \eqref{eq:J0} below, the infected coordinates are listed first. The disease-free set is
\begin{equation}\label{eq:dfset}
\mathcal{X}_0:=\{x\in\Rp^n:\ y=0\}.
\end{equation}
Throughout, $\one$ is the vector of ones and $e_i$ the $i$th standard basis vector of the appropriate dimension, vector and matrix inequalities are componentwise, $u\gneq0$ means $u\ge0$ with $u\ne0$, $\|\cdot\|$ is the Euclidean norm and the induced matrix norm, $\spec(A)$, $\rho(A)$ and $s(A)=\max\{\operatorname{Re}\eta:\eta\in\spec(A)\}$ denote spectrum, spectral radius and spectral abscissa.

\begin{assumption}\label{as:A1}
For each $\theta\in\Theta$ the map $f(\cdot,\theta):\R^n\to\R^n$ is of class $C^2$, and $f$ together with its first two $x$ derivatives is continuous in $(x,\theta)$.
\end{assumption}

\begin{assumption}[Quasi positivity]\label{as:A2}
$f_i(x,\theta)\ge0$ whenever $x\in\Rp^n$, $x_i=0$, $\theta\in\Theta$.
\end{assumption}

\begin{assumption}[Demographic dissipativity]\label{as:A3}
There are constants $\Lambda>0$ and $\mu>0$ such that $\one^{\!\top}f(x,\theta)\le\Lambda-\mu\,\one^{\!\top}x$ for all $x\in\Rp^n$ and $\theta\in\Theta$.
\end{assumption}

Assumption \ref{as:A2} states that no compartment is driven negative from its own boundary face and holds whenever outflows are proportional to compartment contents. It is imposed only on the boundary faces $\{x\in\Rp^n:x_i=0\}$ of the orthant, since only there can a solution leave $\Rp^n$, and in the interior $f_i$ may have either sign; it is required for every $\theta\in\Theta$ so that positivity, and every result below, holds uniformly over the admissible mechanistic parameters. Assumption \ref{as:A3} encodes recruitment at rate at most $\Lambda$ and per capita removal at rate at least $\mu$, with $\mu$ the minimum mortality rate across compartments.

\subsection{Neural closures as fluxes}

The naive universal differential equation $\dot x=f(x,\theta)+\N_\phi(x)$ is unsuitable for epidemic modeling. The following observation applies to any fixed choice of the network parameters $\phi$ for which the stated conditions hold.

\begin{proposition}[Failure of the additive closure]\label{prop:naivefail}
Let $f$ satisfy Assumption \ref{as:A1}, let $\N_\phi:\R^n\to\R^n$ be continuous, put $F:=f(\cdot,\theta)+\N_\phi$, and suppose $F_i(\bar x)<0$ for some $i$ and some $\bar x\in\Rp^n$ with $\bar x_i=0$. Then every solution of $\dot x=F(x)$ with $x(0)=\bar x$ satisfies $x_i(t)<0$ for all sufficiently small $t>0$; in particular $\Rp^n$ is not forward invariant. Solutions exist by Peano's theorem, and they are unique when $\N_\phi$ is locally Lipschitz, for instance a network with a Lipschitz activation (Lemma \ref{lem:nnlip}), but uniqueness is not needed. Moreover, suppose that $\one^{\!\top}f(x,\theta)\ge-\bar\mu\,\one^{\!\top}x$ on $\Rp^n$ for some $\bar\mu>0$, and that $\one^{\!\top}\N_\phi(x)\ge\kappa\|x\|$ for $\|x\|\ge\varrho_1$ with $\kappa>\bar\mu\sqrt n$. Then along every solution that remains in $\Rp^n$ and starts with $\one^{\!\top}x_0\ge\sqrt n\,\varrho_1$, the population $\one^{\!\top}x(t)$ grows at least exponentially.
\end{proposition}

\begin{proof}
Since $F$ is continuous on $\R^n$, there are $r>0$ and $c>0$ with $F_i(x)\le-c$ on $B_r(\bar x)$, and $M\ge c$ with $\|F(x)\|\le M$ there. Let $x$ be a solution with $x(0)=\bar x$ on $[0,T)$, $T>0$, and let $\tau$ be the supremum of the $t\in[0,T)$ with $x([0,t])\subset B_r(\bar x)$. For $t<\tau$ one has $\|x(t)-\bar x\|\le Mt$, so $\tau\ge\min\{T,r/M\}$, since otherwise $x(\tau)$ would lie in the open ball and the solution would remain there beyond $\tau$. Hence $x_i(t)=\int_0^tF_i(x(s))\,\dd s\le-ct<0$ for $0<t<\tau$. For the last claim put $N:=\one^{\!\top}x$. For $x\in\Rp^n$ one has $\|x\|\ge N/\sqrt n$, so $\|x\|\ge\varrho_1$ whenever $N\ge\sqrt n\,\varrho_1$, and then $\dot N=\one^{\!\top}f(x,\theta)+\one^{\!\top}\N_\phi(x)\ge-\bar\mu N+\kappa\|x\|\ge(\kappa/\sqrt n-\bar\mu)N$. The bracket is positive, so $N$ is nondecreasing, the inequality persists, and Gr\"onwall's inequality gives $N(t)\ge N(0)\ee^{(\kappa/\sqrt n-\bar\mu)t}$.
\end{proof}

A compartmental model is determined by a list of fluxes and a stoichiometry: a flux $J_k$ transfers individuals from a source compartment $s_k$ to a target $t_k$, contributing $\nu_kJ_k$ with $\nu_k=e_{t_k}-e_{s_k}$, or $\nu_k=-e_{s_k}$ for removals. Positivity holds because a flux out of compartment $i$ carries the factor $x_i$, and the population is controlled because $\one^{\!\top}\nu_k\le0$. We impose this structure on the learned terms.

\begin{definition}[Admissible closure]\label{def:closure}
Let $q\ge1$ and let $\Phi\subset\R^d$ be compact. An \emph{admissible closure} is a family $\G_\phi=(\G_{\phi,1},\dots,\G_{\phi,q})$ of maps $\G_{\phi,k}:\R^n\to\R$, together with sources $s_k$ and stoichiometric vectors $\nu_k\in\{e_{t_k}-e_{s_k},\,-e_{s_k}\}$, such that each $\G_{\phi,k}$ is of class $C^2$, is continuous in $(x,\phi)$ together with its first two $x$ derivatives, satisfies $\G_{\phi,k}\ge0$ on $\Rp^n$, and is bounded on bounded subsets of $\Rp^n$ uniformly in $\phi\in\Phi$. The closure field is
\begin{equation}\label{eq:closure}
\Xi(x,\phi):=\sum_{k=1}^{q}\nu_k\,x_{s_k}\,\G_{\phi,k}(x).
\end{equation}
\end{definition}

In \eqref{eq:closure}, $\phi\in\Phi$ collects the network weights and biases; it is distinct from the state $x$ and from the mechanistic parameters $\theta$, which are not trained. The outputs of the network are the $q$ scalar per capita rates $\G_{\phi,1}(x),\dots,\G_{\phi,q}(x)$, with units of inverse time: $\G_{\phi,k}(x)$ is the rate at which an individual in the source class $s_k$ undergoes transition $k$, so $x_{s_k}\G_{\phi,k}(x)$ is the number of individuals moved per unit time, and $\nu_k$ records this flux as an outflow from $s_k$ and, unless $k$ is a removal, an inflow to $t_k$. One network with $q$ outputs and $q$ separate networks are both allowed.

A \emph{rate function} $\G_{\phi,k}$ that is globally bounded, $0\le\G_{\phi,k}\le\bar G$ on $\Rp^n$ with $\bar G$ independent of $\phi$, is called \emph{uniformly bounded}. These requirements are architectural rather than restrictive: a feedforward network with a smooth activation is $C^2$ and $C^1$ in its weights, and composing the output with a scaled logistic map $u\mapsto\bar G(1+\ee^{-u})^{-1}$ enforces uniform boundedness exactly, at the cost of an a priori ceiling $\bar G$ on the learned rates, which is available in epidemiological practice. Boundedness on bounded sets, rather than global boundedness, is required in Definition \ref{def:closure} because the transmission rate introduced next contains the infected coordinates multiplicatively and need not be globally bounded.

\subsection{Epidemiological consistency}

If the closure injects individuals into infected compartments when none are present, the set $\mathcal X_0$ ceases to be invariant, no disease-free equilibrium exists, and no reproduction number can be defined. Let $\iota_{\Ical}:\R^m\to\R^n$ denote the embedding of the infected coordinates.

\begin{assumption}[Epidemiological consistency]\label{as:A4}
$f_i(z,0,\theta)=0$ for all $i\in\Ical$, $z\in\Rp^{n-m}$, $\theta\in\Theta$. The closure splits as $\Xi=\Xi^{\mathrm{tr}}+\Xi^{\mathrm{pr}}+\Xi^{\mathrm{un}}$, where the transmission part is
\begin{equation}\label{eq:foi}
\Xi^{\mathrm{tr}}(x,\phi)=\bigl(\hat w-e_{s_\ast}\bigr)\,x_{s_\ast}\,\lambda_\phi(x),\qquad
\lambda_\phi(x)=\sum_{j\in\Ical}\beta^{\phi}_j(x)\,x_j,\qquad \hat w:=\iota_{\Ical}w,
\end{equation}
with $s_\ast\in\Ucal$ the susceptible index, $w\in\Rp^{m}$ a fixed distribution vector with $\one^{\!\top}w=1$ describing the classes entered upon infection, and $\beta^\phi=(\beta^\phi_j)_{j\in\Ical}$ a family of uniformly bounded rate functions, $0\le\beta^\phi_j\le\bar G$ on $\Rp^n$, satisfying the regularity requirements of Definition \ref{def:closure}; the progression part $\Xi^{\mathrm{pr}}$ consists of fluxes with source in $\Ical$; and the uninfected part $\Xi^{\mathrm{un}}$ consists of fluxes with source in $\Ucal$ whose target, if any, lies in $\Ucal$.
\end{assumption}

The transmission flux is the sum over $j\in\Ical$ of the admissible fluxes $(e_j-e_{s_\ast})\,x_{s_\ast}\,w_j\lambda_\phi(x)$, whose rates $w_j\lambda_\phi$ are nonnegative, $C^2$ and bounded on bounded sets but not globally bounded; this is why Definition \ref{def:closure} only requires boundedness on bounded sets. What the factorized form \eqref{eq:foi} buys is that $\lambda_\phi(z,0)=0$ identically in $\phi$. Mass action $\beta SI$, saturating incidence $\beta SI/(1+\alpha I)$ \cite{CapassoSerio1978} and nonmonotone incidence $\beta SI/(1+\alpha I^2)$ \cite{XiaoRuan2007} are of this form with $\beta_I$ equal to $\beta$, $\beta/(1+\alpha I)$ and $\beta/(1+\alpha I^2)$ respectively. More generally, sufficiently regular forces of infection that vanish on $\mathcal X_0$ and are nondecreasing in the infected compartments, with partial derivatives in those compartments bounded by $\bar G$, admit the representation (Proposition \ref{prop:factor}).

\begin{proposition}[Factorization through the infected coordinates]\label{prop:factor}
Let $\lambda:\R^n\to\R$ be of class $C^3$ with $\lambda(z,0)=0$ for all $z$. Then $\lambda(z,y)=\sum_{j\in\Ical}\beta_j(z,y)\,y_j$ with
\begin{equation}\label{eq:hadamard}
\beta_j(z,y):=\int_0^1\partial_{y_j}\lambda(z,\tau y)\,\dd\tau ,
\end{equation}
and each $\beta_j$ is of class $C^2$. If moreover $0\le\partial_{y_j}\lambda\le\bar G$ on $\Rp^n$ for each $j\in\Ical$, then $0\le\beta_j\le\bar G$ on $\Rp^n$, so $(\beta_j)_{j\in\Ical}$ satisfies the requirements placed on $\beta^\phi$ in Assumption \ref{as:A4}. If $\lambda=\lambda_\phi$ depends on $\phi$ with the derivatives $\partial_{y_j}\lambda_\phi$ and their first two $x$ derivatives continuous in $(x,\phi)$, the same holds for $\beta_j$.
\end{proposition}

\begin{proof}
For fixed $(z,y)$ the map $\tau\mapsto\lambda(z,\tau y)$ is $C^1$ on $[0,1]$ with derivative $\sum_j\partial_{y_j}\lambda(z,\tau y)y_j$, and it vanishes at $\tau=0$; integrating gives the representation. Differentiation under the integral sign, justified by the continuity of the third derivatives of $\lambda$ on the compact segment, gives $\beta_j\in C^2$ and the continuity in $(x,\phi)$. The bounds follow from those on the integrand, since the segment $\{(z,\tau y):\tau\in[0,1]\}$ lies in $\Rp^n$ when $(z,y)$ does.
\end{proof}

The representation \eqref{eq:hadamard} is not unique, and Assumption \ref{as:A4} does not require monotonicity: the nonmonotone incidence above is admissible directly, with a coefficient $\beta_I$ that decreases in $I$. Forces of infection whose partial derivatives with respect to the infected compartments change sign, or which are less regular than $C^3$, are covered by \eqref{eq:foi} whenever a nonnegative $C^2$ factorization can be written down, and not otherwise. The uninfected part allows learned processes among uninfected classes, such as waning immunity from a recovered to the susceptible class; it is absent in the running examples below.

Assumption \ref{as:A4} restricts the transmission architecture in two ways: all new infections come from the single susceptible class $s_\ast$, and they are distributed among the infected classes by a fixed vector $w$ that depends neither on the state nor on the source of infection. Models with several susceptible classes, such as vaccinated and unvaccinated groups or multiple host types, vector-host models, in which susceptible hosts and susceptible vectors are infected through different routes, and models in which the class entered upon infection depends on the route are therefore outside the class; Remark \ref{rem:rankone} describes the consequence for the next-generation matrix.

Collecting the pieces, the object of study is
\begin{equation}\label{eq:ude}
\dot x=F_\phi(x):=f(x,\theta)+\sum_{k=1}^{q}\nu_k\,x_{s_k}\,\G_{\phi,k}(x),\qquad x(0)=x_0\in\Rp^n,
\end{equation}
under Assumptions \ref{as:A1}--\ref{as:A4}. We call \eqref{eq:ude} an \emph{epidemic UDE} and denote the class of such systems by $\Ecal$. Not every result below needs every assumption; Table \ref{tab:assumptions} records, for each group of results, the hypotheses actually used, including the assumptions on the disease-free equilibrium (Assumption \ref{as:A5}), the transfer matrix (Assumption \ref{as:V}) and the endemic branch (Assumption \ref{as:A6}) introduced later, and the conditions on the closure that are imposed or certified after training.

\begin{table}[t]
\centering
\caption{Hypotheses used by each group of results. D denotes Definition \ref{def:closure}; A1--A6 and V the assumptions of that name; Dom, Sub, Mon, Lip the domination, sublinearity, monotonicity and Lipschitz conditions on the closure.}
\label{tab:assumptions}
\footnotesize
\begin{tabular}{@{}lll@{}}
\toprule
Result & Location & Hypotheses used \\
\midrule
Existence, uniqueness, dependence on data & Thms \ref{thm:local}, \ref{thm:contdep} & A1, D \\
Positivity & Thm \ref{thm:pos} & A1, A2, D \\
Dissipativity, feasible region & Thm \ref{thm:region} & A1--A3, D \\
Invariance of the disease-free set & Prop \ref{prop:dfinv} & A1, A2, A4, D \\
Next-generation matrix, closed form $\Rz$, bounds & Prop \ref{prop:ngm}, Thms \ref{thm:R0closed}, \ref{thm:R0bounds} & A1--A5, V \\
Local stability of $E_0$ & Thm \ref{thm:LASDFE} & A1--A5, V \\
Global stability of $E_0$, extinction & Thm \ref{thm:GASDFE}, Cor \ref{cor:extinct} & A1--A5, V, Dom \\
Endemic equilibrium: existence & Thm \ref{thm:EEexist} & A1--A6, V \\
Endemic equilibrium: uniqueness & Thm \ref{thm:EEunique} & A1--A6, V, Sub \\
Endemic equilibrium: global stability & Thm \ref{thm:GASEE} & class $\Ecal_{\mathrm{sep}}$, Mon, Sub \\
Uniform persistence & Thm \ref{thm:persist} & A1--A5, V, $F-V$ irreducible \\
Bifurcation at threshold & Thm \ref{thm:backward} & model \eqref{eq:seir}, $\beta_\phi\in C^2$ \\
Sensitivity & Sec \ref{sec:sensitivity} & as Thm \ref{thm:R0closed} \\
Nonidentifiability of $\Rz$ and of the bifurcation type & Thm \ref{thm:nonid}, Cor \ref{cor:bifnonid} & A1--A5, V \\
Quantitative and monotone identifiability bounds & Thm \ref{thm:quantid}, Prop \ref{prop:monoident} & A1--A5, V, Lip or Mon \\
\bottomrule
\end{tabular}
\end{table}

\begin{example}[Neural SEIR model]\label{ex:seir}
Take $n=4$, $x=(S,E,I,R)$, $\Ical=\{2,3\}$, $s_\ast=1$, $w=(1,0)^{\!\top}$, and let the closure supply $\lambda_\phi(x)=\beta_\phi(x)I$ with $0\le\beta_\phi\le\bar G$:
\begin{equation}\label{eq:seir}
\dot S=\Lambda-\beta_\phi(x)SI-\mu S,\quad
\dot E=\beta_\phi(x)SI-(\mu+\sigma)E,\quad
\dot I=\sigma E-(\mu+\gamma+\delta)I,\quad
\dot R=\gamma I-\mu R,
\end{equation}
with $\Lambda,\mu,\sigma,\gamma>0$ and $\delta\ge0$. The learned quantity is the scalar rate $\beta_\phi(x)$, so $\beta^\phi=(\beta^\phi_E,\beta^\phi_I)=(0,\beta_\phi)$; the closure field it induces is the vector $\Xi^{\mathrm{tr}}(x,\phi)=(-1,1,0,0)^{\!\top}S\beta_\phi(x)I$, which enters the $S$ and $E$ equations with opposite signs, and all other terms of \eqref{eq:seir} are mechanistic. Assumptions \ref{as:A1}--\ref{as:A4} hold, since $\one^{\!\top}f=\Lambda-\mu\one^{\!\top}x-\delta I\le\Lambda-\mu\one^{\!\top}x$. This model is the running illustration, and the computations of Sections \ref{sec:bifurcation} and \ref{sec:sensitivity} are carried out for it.
\end{example}

\begin{example}[Separable subclass]\label{ex:sir}
Take $n=3$, $x=(S,I,R)$, $\Ical=\{2\}$, and $\lambda_\phi(x)=h_\phi(I)$ with $h_\phi(I)=I\beta_\phi(I)$:
\begin{equation}\label{eq:sir}
\dot S=\Lambda-Sh_\phi(I)-\mu S,\qquad \dot I=Sh_\phi(I)-(\mu+\gamma+\delta)I,\qquad \dot R=\gamma I-\mu R.
\end{equation}
We write $\Ecal_{\mathrm{sep}}$ for this subclass; it is the setting of Theorem \ref{thm:GASEE}.
\end{example}

\section{Well-posedness, positivity, and invariance}\label{sec:wp}

\subsection{Regularity}

\begin{lemma}\label{lem:nnlip}
Let $\N(x)=W_L\varsigma(W_{L-1}\varsigma(\cdots\varsigma(W_1x+b_1)\cdots)+b_{L-1})+b_L$ with $\varsigma$ acting componentwise and Lipschitz with constant $L_\varsigma$. Then $\N$ is globally Lipschitz with constant $L_\varsigma^{L-1}\prod_{\ell=1}^{L}\|W_\ell\|$ and grows at most linearly; it is of class $C^r$ if $\varsigma\in C^r$; and if the output is composed with a bounded Lipschitz $\varrho$, then $\varrho\circ\N$ is bounded and globally Lipschitz.
\end{lemma}

\begin{proof}
Since $\varsigma$ acts componentwise, $\|\varsigma(a)-\varsigma(\tilde a)\|^2=\sum_i|\varsigma(a_i)-\varsigma(\tilde a_i)|^2\le L_\varsigma^2\|a-\tilde a\|^2$. Writing $u^{(0)}(x)=x$ and $u^{(\ell)}(x)=\varsigma(W_\ell u^{(\ell-1)}(x)+b_\ell)$, induction gives $\|u^{(L-1)}(x)-u^{(L-1)}(\tilde x)\|\le L_\varsigma^{L-1}\prod_{\ell<L}\|W_\ell\|\,\|x-\tilde x\|$, and applying $W_L$ gives the constant. Linear growth follows from $\|\N(x)\|\le\|\N(0)\|+\mathrm{Lip}(\N)\|x\|$, smoothness from the chain rule, and the last claim from composition of Lipschitz maps.
\end{proof}

\begin{corollary}\label{cor:xireg}
Let the rate functions of Definition \ref{def:closure} be of the form $\varrho\circ\N_{\phi,k}$ with $\varrho$ bounded, $C^2$ and Lipschitz, or of the form $w_j\lambda_\phi$ with $\beta^\phi_j=\varrho\circ\N_{\phi,j}$ as in \eqref{eq:foi}. Then for every $\varrho_0>0$ there are constants $G_{\varrho_0},L_{\varrho_0}>0$, independent of $\phi\in\Phi$, such that $0\le\G_{\phi,k}\le G_{\varrho_0}$ and $\mathrm{Lip}(\G_{\phi,k})\le L_{\varrho_0}$ on $\{x\in\Rp^n:\|x\|\le\varrho_0\}$. Moreover, $\Xi(\cdot,\phi)$, which is $C^2$ by Definition \ref{def:closure}, satisfies $\|\Xi(x,\phi)-\Xi(\tilde x,\phi)\|\le\sqrt2\,q\,(G_{\varrho_0}+\varrho_0L_{\varrho_0})\|x-\tilde x\|$ whenever $x,\tilde x\in\Rp^n$ and $\|x\|,\|\tilde x\|\le\varrho_0$.
\end{corollary}

\begin{proof}
For rates of the form $\varrho\circ\N_{\phi,k}$, global bounds $\bar G$ and $L_G$ independent of $\phi$ follow from Lemma \ref{lem:nnlip}, continuity of the weights in $\phi$ and compactness of $\Phi$. For $w_j\lambda_\phi=w_j\sum_i\beta^\phi_i x_i$ one has $0\le w_j\lambda_\phi\le\bar G\sqrt m\,\varrho_0$ and $\mathrm{Lip}\le\bar G\sqrt m+\varrho_0\sqrt m\,L_G$ on the ball, again uniformly in $\phi$. For the last bound, $|x_{s_k}\G_{\phi,k}(x)-\tilde x_{s_k}\G_{\phi,k}(\tilde x)|\le G_{\varrho_0}|x_{s_k}-\tilde x_{s_k}|+\varrho_0L_{\varrho_0}\|x-\tilde x\|$, and $\|\nu_k\|\le\sqrt2$; summing over the $q$ fluxes gives the claim.
\end{proof}

\subsection{Existence, uniqueness and dependence on data}

\begin{theorem}\label{thm:local}
Under Assumption \ref{as:A1} and Definition \ref{def:closure}, the vector field $F_\phi$ of \eqref{eq:ude} is of class $C^2$ on $\R^n$, and for every $x_0\in\R^n$ the equation $\dot x=F_\phi(x)$ with $x(0)=x_0$ has a unique solution on a maximal interval $(T_-(x_0),T_+(x_0))$; moreover $\|x(t)\|\to\infty$ as $t\uparrow T_+(x_0)$ if $T_+(x_0)<\infty$. In particular this holds for the initial data $x_0\in\Rp^n$ of \eqref{eq:ude}.
\end{theorem}

\begin{proof}
By Assumption \ref{as:A1}, $f(\cdot,\theta)$ is $C^2$ on $\R^n$, and by Definition \ref{def:closure} every rate $\G_{\phi,k}$, including the transmission rates of Assumption \ref{as:A4}, is $C^2$ on $\R^n$. Hence $F_\phi\in C^2(\R^n;\R^n)$, and in particular $F_\phi$ is locally Lipschitz. The sign and boundedness requirements of Definition \ref{def:closure}, which concern $\Rp^n$, are not used here; they enter in Theorems \ref{thm:pos} and \ref{thm:region}. Fix $x_0$, let $B:=\bar B_1(x_0)$, let $L_B>0$ and $M_B>0$ be upper bounds for the Lipschitz constant and the norm of $F_\phi$ on $B$, and set $a:=\min\{1/(2L_B),1/(2M_B)\}$. On the complete metric space $X$ of continuous maps $[-a,a]\to B$ with the supremum metric, the Picard operator $(\mathcal{T}u)(t)=x_0+\int_0^tF_\phi(u(s))\dd s$ satisfies $\|(\mathcal Tu)(t)-x_0\|\le aM_B\le1/2$, so it maps $X$ into $X$, and $d(\mathcal Tu,\mathcal Tv)\le aL_Bd(u,v)\le d(u,v)/2$, so it is a contraction. Its unique fixed point is the unique solution with values in $B$; uniqueness among all solutions follows from Gr\"onwall's inequality applied on compact subintervals, where both solutions are bounded and $F_\phi$ is Lipschitz. Maximal continuation is standard for locally Lipschitz fields \cite[Ch.~II]{Hartman2002}.
\end{proof}

\begin{theorem}\label{thm:contdep}
Let $x(\cdot;x_0,\phi)$ and $x(\cdot;\tilde x_0,\tilde\phi)$ solve \eqref{eq:ude} and remain in a compact set $K$ on $[0,T]$, and let $L_K>0$ be a Lipschitz constant for $F_\phi$ on $K$. Then for $t\in[0,T]$,
\begin{equation}\label{eq:contdep}
\|x(t;x_0,\phi)-x(t;\tilde x_0,\tilde\phi)\|\le \ee^{L_Kt}\|x_0-\tilde x_0\|+\frac{\ee^{L_Kt}-1}{L_K}\sup_{x\in K}\|F_\phi(x)-F_{\tilde\phi}(x)\|.
\end{equation}
If moreover $\phi\mapsto\G_{\phi,k}(x)$ is $C^1$ with derivative continuous in $(x,\phi)$, then $\phi\mapsto x(t;x_0,\phi)$ is $C^1$ and $\Psi(t):=D_\phi x(t;x_0,\phi)$ solves $\dot\Psi=D_xF_\phi(x(t))\Psi+D_\phi F_\phi(x(t))$, $\Psi(0)=0$.
\end{theorem}

\begin{proof}
Set $\psi(t):=\|x(t;x_0,\phi)-x(t;\tilde x_0,\tilde\phi)\|$ and $\varepsilon:=\sup_K\|F_\phi-F_{\tilde\phi}\|$. The integral formulation gives $\psi(t)\le c+\int_0^t(\varepsilon+L_K\psi(s))\dd s=:\chi(t)$ with $c=\|x_0-\tilde x_0\|$. Then $\chi'=\varepsilon+L_K\psi\le\varepsilon+L_K\chi$, so $(\ee^{-L_Kt}\chi)'\le\varepsilon\ee^{-L_Kt}$, and integrating yields $\chi(t)\le c\ee^{L_Kt}+\varepsilon(\ee^{L_Kt}-1)/L_K$, which is \eqref{eq:contdep} since $\psi\le\chi$. The second statement is the classical theorem on differentiable dependence on parameters for fields that are $C^1$ jointly in $(x,\phi)$ \cite[Ch.~V]{Hartman2002}, the derivative being characterized by the variational equation, which is linear with continuous coefficients and hence uniquely solvable on the same interval.
\end{proof}

The second part of Theorem \ref{thm:contdep} legitimizes gradient based training: the sensitivities computed by adjoint solvers \cite{Chen2018NeuralODE,RackauckasNie2017} are genuine derivatives of the flow.

\subsection{Positivity, boundedness, and the feasible region}

\begin{theorem}[Positivity]\label{thm:pos}
Under Assumptions \ref{as:A1}, \ref{as:A2} and Definition \ref{def:closure}, $\Rp^n$ is forward invariant for \eqref{eq:ude} for every $\phi\in\Phi$, and solutions starting in the interior remain in the interior.
\end{theorem}

\begin{proof}
We first show that $F_\phi$ is quasi positive on $\Rp^n$. Let $x\in\Rp^n$ with $x_i=0$. Then $f_i(x,\theta)\ge0$ by Assumption \ref{as:A2}. In \eqref{eq:closure}, the $i$th component of $\nu_kx_{s_k}\G_{\phi,k}(x)$ is nonnegative unless $s_k=i$, because $(\nu_k)_i\ge0$ in that case while $x_{s_k}\ge0$ and $\G_{\phi,k}\ge0$; and if $s_k=i$ the term equals $-x_i\G_{\phi,k}(x)=0$. Hence $F_{\phi,i}(x)\ge0$.

Let $x_0\in\Rp^n$ and $[0,T]\subset[0,T_+(x_0))$, and let $K$ be a compact neighborhood of the trajectory on $[0,T]$. For $\epsilon>0$ consider $\dot x^\epsilon=F_\phi(x^\epsilon)+\epsilon\one$, $x^\epsilon(0)=x_0+\epsilon\one$. The Gr\"onwall argument in the proof of Theorem \ref{thm:contdep} applies verbatim to the pair of locally Lipschitz fields $F_\phi$ and $F_\phi+\epsilon\one$, whose difference has norm $\epsilon\sqrt n$; hence there is $\epsilon_0>0$ such that for $\epsilon\in(0,\epsilon_0]$ the solution $x^\epsilon$ exists on $[0,T]$, stays in $K$ and satisfies $\sup_{[0,T]}\|x^\epsilon-x\|\le C\epsilon$. If $x^\epsilon$ left the open orthant, then $t_1:=\inf\{t\in[0,T]:\min_ix_i^\epsilon(t)\le0\}$ would satisfy $t_1>0$, $x^\epsilon(t_1)\ge0$ and $x^\epsilon_i(t_1)=0$ for some $i$; quasi positivity gives $\dot x_i^\epsilon(t_1)=F_{\phi,i}(x^\epsilon(t_1))+\epsilon\ge\epsilon>0$, whereas $x_i^\epsilon(t)>0=x^\epsilon_i(t_1)$ for $t<t_1$ forces $\dot x^\epsilon_i(t_1)\le0$, a contradiction. Hence $x^\epsilon>0$ on $[0,T]$, and letting $\epsilon\downarrow0$ gives $x\ge0$ on $[0,T]$; as $T<T_+(x_0)$ was arbitrary, $\Rp^n$ is invariant.

If $x_0>0$, suppose $\tau:=\inf\{t\ge0:x_i(t)=0\ \text{for some}\ i\}<T_+(x_0)$. On $[0,\tau]$ the solution lies in a compact box $K':=\{x\in\Rp^n:x\le b\}$ for some $b>0$, on which $F_\phi$ is Lipschitz with constant $L'$. For $x\in K'$ let $\hat x$ be $x$ with the $i$th coordinate set to zero; then $\hat x\in K'$, $\|x-\hat x\|=x_i$ and $F_{\phi,i}(\hat x)\ge0$, so $F_{\phi,i}(x)\ge-L'x_i$. Hence $\dot x_i\ge-L'x_i$ and $x_i(t)\ge x_i(0)\ee^{-L't}>0$ on $[0,\tau]$, contradicting $x_i(\tau)=0$.
\end{proof}

\begin{theorem}[Dissipativity and feasible region]\label{thm:region}
Let Assumptions \ref{as:A1}--\ref{as:A3} and Definition \ref{def:closure} hold, let $x_0\in\Rp^n$ and $N(t):=\one^{\!\top}x(t)$. Then $T_+(x_0)=\infty$ and
\begin{equation}\label{eq:Nbound}
N(t)\le N(0)\ee^{-\mu t}+\frac{\Lambda}{\mu}\bigl(1-\ee^{-\mu t}\bigr),\qquad t\ge0 .
\end{equation}
The set
\begin{equation}\label{eq:Omega}
\Om:=\Bigl\{x\in\Rp^n:\ \one^{\!\top}x\le\Lambda/\mu\Bigr\}
\end{equation}
is compact, convex, forward invariant for every $\phi\in\Phi$, and attracts every bounded subset of $\Rp^n$; for each $\eta>0$ the compact set $\mathcal A_\eta:=\{x\in\Rp^n:\one^{\!\top}x\le\Lambda/\mu+\eta\}$ is forward invariant and absorbing. Every compact invariant subset of $\Rp^n$, in particular every equilibrium, lies in $\Om$, and $0\le x_i(t)\le\Lambda/\mu$ for all $i$ when $x_0\in\Om$.
\end{theorem}

\begin{proof}
By Theorem \ref{thm:pos} the solution is nonnegative. Summing \eqref{eq:ude} and using $\one^{\!\top}\nu_k\le0$ with $x_{s_k}\G_{\phi,k}(x)\ge0$ gives $\dot N\le\one^{\!\top}f(x,\theta)\le\Lambda-\mu N$ by Assumption \ref{as:A3}; multiplying by $\ee^{\mu t}$ and integrating gives \eqref{eq:Nbound}. Hence $\|x(t)\|\le N(t)\le\max\{N(0),\Lambda/\mu\}$, so the solution is bounded and $T_+(x_0)=\infty$ by Theorem \ref{thm:local}. Invariance of $\Om$ and of $\mathcal A_\eta$ is immediate from \eqref{eq:Nbound}, and if $N_B$ bounds $\one^{\!\top}x_0$ on a bounded set $B$ then $N(t)\le\Lambda/\mu+N_B\ee^{-\mu t}$, giving absorption in time $\mu^{-1}\log^{+}(N_B/\eta)$. For attraction, if $x\ge0$ and $N=\one^{\!\top}x>\Lambda/\mu$ then $\tilde x:=(\Lambda/(\mu N))x\in\Om$ and $\|x-\tilde x\|=\|x\|(1-\Lambda/(\mu N))\le N-\Lambda/\mu$, so $\dist(x(t),\Om)\le N(0)\ee^{-\mu t}$ by \eqref{eq:Nbound}. If $\mathcal K\subset\Rp^n$ is compact and invariant and $\bar x\in\mathcal K$, then for each $t>0$ there is $x_{-t}\in\mathcal K$ with $\bar x=x(t;x_{-t})$, and \eqref{eq:Nbound} gives $\one^{\!\top}\bar x\le(\sup_{\mathcal K}\one^{\!\top}x)\ee^{-\mu t}+(\Lambda/\mu)(1-\ee^{-\mu t})\to\Lambda/\mu$. Finally $x_i\le\one^{\!\top}x$ for $x\ge0$.
\end{proof}

Theorem \ref{thm:region} shows that the feasible region of the mechanistic model remains feasible for the whole learned family, and identifies $\Om$ as the compact set on which all uniform conditions on the closure are to be imposed and verified.

\section{Disease-free equilibrium and the reproduction number}\label{sec:R0}

\begin{proposition}\label{prop:dfinv}
Under Assumptions \ref{as:A1}, \ref{as:A2}, \ref{as:A4} and Definition \ref{def:closure} the set $\mathcal X_0$ is forward invariant for every $\phi\in\Phi$, and on $\mathcal X_0$ the dynamics reduce to
\begin{equation}\label{eq:dfsub}
\dot z=f^{\Ucal}(z,0,\theta)+\Xi^{\Ucal}((z,0),\phi)=:g_\phi(z),
\end{equation}
where the superscript denotes restriction to the uninfected coordinates.
\end{proposition}

\begin{proof}
Let $i\in\Ical$. Then $f_i(z,0,\theta)=0$ by Assumption \ref{as:A4}; transmission fluxes carry the factor $\lambda_\phi(z,0)=0$; progression fluxes have source in $\Ical$ and carry a factor $x_{s_k}=0$; and fluxes of the uninfected part have source and target in $\Ucal$, so their $i$th component vanishes identically. Hence $F_{\phi,i}(z,0)=0$ for $i\in\Ical$, so the field is tangent to $\mathcal X_0$. The reduced field $g_\phi$ is quasi positive on $\Rp^{n-m}$, because its components are those of $F_\phi(z,0)$ in the uninfected coordinates and $F_\phi$ is quasi positive on $\Rp^n$ (proof of Theorem \ref{thm:pos}); the argument of Theorem \ref{thm:pos} therefore shows that solutions of \eqref{eq:dfsub} with $z(0)\ge0$ stay nonnegative. For such a solution, $(z(t),0)$ solves \eqref{eq:ude}, and uniqueness of solutions (Theorem \ref{thm:local}) identifies it with the solution from $(z(0),0)$. Hence $\mathcal X_0$ is forward invariant, and the dynamics on it are given by \eqref{eq:dfsub}.
\end{proof}

\begin{assumption}\label{as:A5}
The reduced system \eqref{eq:dfsub} has a unique equilibrium $z^\ast=z^\ast(\theta,\phi)\ge0$ in the projection of $\Om$ onto the uninfected coordinates, with a positive susceptible component $z^\ast_{s_\ast}>0$, and $z^\ast$ is globally asymptotically stable there for \eqref{eq:dfsub}, with $s(Dg_\phi(z^\ast))<0$.
\end{assumption}

We call $E_0:=(z^\ast,0)\in\Om$ the disease-free equilibrium. Only the positivity of the susceptible component $z^\ast_{s_\ast}$ is used below; other uninfected classes, such as a recovered class, are typically empty at $E_0$.

\begin{proposition}\label{prop:dfeindep}
Let $f^{\Ucal}$ denote the restriction of the mechanistic vector field $f$ to the uninfected coordinates. If the uninfected part $\Xi^{\mathrm{un}}$ of the closure vanishes, then $g_\phi=f^{\Ucal}(\cdot,0,\theta)$ is independent of $\phi$, and so is $E_0$. For \eqref{eq:seir}, $E_0=(\Lambda/\mu,0,0,0)$ for every $\phi$, and Assumption \ref{as:A5} holds.
\end{proposition}

\begin{proof}
On $\mathcal X_0$ transmission fluxes vanish by \eqref{eq:foi} and progression fluxes vanish because their sources lie in $\Ical$; if $\Xi^{\mathrm{un}}=0$ no other closure flux remains. For \eqref{eq:seir} the reduced system on $\mathcal X_0$ is
\[
\dot S=\Lambda-\mu S,\qquad \dot R=-\mu R.
\]
It has the unique globally stable equilibrium $(\Lambda/\mu,0)$, whose susceptible component is positive. Its Jacobian with respect to the uninfected variables $(S,R)$ is
\[
Dg_\phi(S,R)=
\begin{pmatrix}
-\mu&0\\
0&-\mu
\end{pmatrix}
=-\mu I_2,
\]
so Assumption \ref{as:A5} holds.
\end{proof}

If instead the closure contains, say, a learned waning immunity flux from the recovered to the susceptible class, then $E_0$ depends on $\phi$ and must be recomputed whenever the network is retrained.

Let $F_\phi^{\Ical}(x)\in\R^m$ denote the restriction of the full vector field $F_\phi(x)$ to the infected coordinates. Write the infected subsystem of \eqref{eq:ude} as
\[
\dot y=\Fcal_\phi(x)-\Vcal_\phi(x),
\]
where $\Fcal_\phi,\Vcal_\phi:\Rp^n\to\R^m$ are defined by
\begin{equation}\label{eq:FV}
\Fcal_\phi(x):=w\,x_{s_\ast}\lambda_\phi(x),\qquad
\Vcal_\phi(x):=\Fcal_\phi(x)-F^{\Ical}_\phi(x).
\end{equation}
Thus $\Fcal_\phi$ collects new infections, whereas $\Vcal_\phi$ collects the net effect of all other transfers in the infected subsystem. We then set
\begin{equation}\label{eq:FVmat}
F(\phi):=D_y\Fcal_\phi(E_0),\qquad V(\phi):=D_y\Vcal_\phi(E_0)\ \in\R^{m\times m}.
\end{equation}
Here $D_y$ is the partial Jacobian with respect to the infected coordinates $y$ at fixed $z$, so $F_{ij}(\phi)=\partial\Fcal_{\phi,i}/\partial y_j(E_0)$ and $V_{ij}(\phi)=\partial\Vcal_{\phi,i}/\partial y_j(E_0)$ for $i,j\in\Ical$. The matrices $F(\phi)$ and $V(\phi)$ are $m\times m$ because they are the Jacobians of the $m$ infected equations with respect to the $m$ infected coordinates $y$.

\begin{assumption}[Transfer matrix]\label{as:V}
$s(-V(\phi))<0$ for every $\phi\in\Phi$.
\end{assumption}

Assumption \ref{as:V} states that, in the absence of new infections, the infected classes are emptied; it is the standing hypothesis (A5) of \cite{vandenDriesscheWatmough2002}. Since $V(\phi)$ has nonpositive off diagonal entries (Proposition \ref{prop:ngm} below), it holds whenever $V(\phi)$ is strictly column diagonally dominant, $\one^{\!\top}V(\phi)>0$ componentwise, which is the case when every infected class is subject to positive removal and all transfers among infected classes are conservative \cite[Ch.~6]{BermanPlemmons1994}. For \eqref{eq:seir} it holds because $V$ is lower triangular with positive diagonal, see Example \ref{ex:R0seir}.

\begin{proposition}\label{prop:ngm}
Under Assumptions \ref{as:A1}--\ref{as:A5},
\begin{equation}\label{eq:Frank1}
F(\phi)=z^\ast_{s_\ast}\,w\,\bigl(\beta^\phi(E_0)\bigr)^{\!\top}\ \ge0
\end{equation}
has rank at most one, and $V(\phi)$ has nonpositive off diagonal entries. If in addition Assumption \ref{as:V} holds, then $V(\phi)$ is a nonsingular M matrix, $V(\phi)^{-1}\ge0$, and $V(\phi)^{-1}$ has positive diagonal.
\end{proposition}

\begin{proof}
By \eqref{eq:foi}, $\Fcal_\phi(x)=w\,x_{s_\ast}\sum_j\beta_j^\phi(x)x_j$; differentiating in $y$ at $E_0$, where $y=0$, annihilates the terms containing $D_y\beta^\phi$ and leaves \eqref{eq:Frank1}, which is nonnegative because $z^\ast_{s_\ast}>0$, $w\ge0$, $\beta^\phi\ge0$. For $i\ne j$ in $\Ical$, \eqref{eq:FV} gives $V_{ij}=F_{ij}-\partial_{y_j}F^{\Ical}_{\phi,i}(E_0)$, and we show that $\partial_{y_j}F^{\Ical}_{\phi,i}(E_0)\ge F_{ij}$. The mechanistic contribution is $\partial_{y_j}f_i(E_0)$: the map $t\mapsto f_i(z^\ast,te_j,\theta)$, $t\ge0$, is nonnegative by Assumption \ref{as:A2}, since the $i$th coordinate of $(z^\ast,te_j)$ vanishes, and it vanishes at $t=0$ by Assumption \ref{as:A4}; hence its derivative at $t=0$ is nonnegative. Among the closure fluxes, those with source $i$ contribute $-x_i\G_{\phi,k}(x)$, whose $y_j$ derivative at $E_0$ vanishes because $x_i=0$ there; those with target $i$ and source in $\Ical$ contribute $x_{s_k}\G_{\phi,k}(x)$, whose $y_j$ derivative at $E_0$ equals $\delta_{s_kj}\G_{\phi,k}(E_0)\ge0$ because $x_{s_k}=0$ there; and the transmission flux contributes $w_ix_{s_\ast}\lambda_\phi(x)$, whose $y_j$ derivative at $E_0$ is exactly $F_{ij}$. No other flux has a nonzero $i$th component. Summing, $V_{ij}\le0$. A matrix with nonpositive off diagonal entries and $s(-V)<0$ is a nonsingular M matrix with nonnegative inverse \cite[Ch.~6]{BermanPlemmons1994}; writing $V=s_0I-B$ with $B\ge0$ and $s_0>\rho(B)$ gives $V^{-1}=s_0^{-1}\sum_{k\ge0}(B/s_0)^k\ge s_0^{-1}I$, so the diagonal of $V^{-1}$ is positive.
\end{proof}
\begin{remark}[Scope of the rank-one structure]\label{rem:rankone}
The rank-one form in \eqref{eq:Frank1} is a consequence of the single-entry transmission architecture in Assumption \ref{as:A4}: infections arise from one susceptible source $s_\ast$ and enter the infected subsystem according to the fixed distribution vector $w$. It is therefore a structural property of the epidemic UDE class considered here, not of arbitrary compartmental epidemic models. Models with multiple susceptible groups, spatial patches, host types, or distinct infection-entry mechanisms generally lead to a higher-rank new-infection matrix and require the corresponding higher-rank next-generation formulation. For example, in a host-vector model infected vectors generate new host infections and infected hosts generate new vector infections; the new-infection matrix then has two nonzero off-diagonal blocks and rank two in general, and $\Rz$ is the spectral radius of a next-generation matrix with zero diagonal blocks, in the simplest case the geometric mean of the host-to-vector and vector-to-host reproduction numbers, rather than \eqref{eq:R0closed}. The results below that rely on \eqref{eq:Frank1} or \eqref{eq:R0closed} are not claimed for such models.
\end{remark}
\begin{definition}\label{def:R0}
The next-generation matrix of \eqref{eq:ude} is $K(\phi):=F(\phi)V(\phi)^{-1}$ and the basic reproduction number is $\Rz(\phi):=\rho(K(\phi))$.
\end{definition}

\begin{theorem}[Closed form and locality]\label{thm:R0closed}
For every epidemic UDE in the single-entry admissible class defined above (Remark \ref{rem:rankone}) that satisfies Assumptions \ref{as:A1}--\ref{as:A5} and \ref{as:V},

\begin{equation}\label{eq:R0closed}
\Rz(\phi)=z^{\ast}_{s_\ast}\,\bigl(\beta^{\phi}(E_0)\bigr)^{\!\top}V(\phi)^{-1}w .
\end{equation}
For fixed $z^\ast$ and $V(\phi)$, $\Rz$ is therefore a nonnegative linear functional of the transmission value $\beta^\phi(E_0)\in\Rp^m$. More generally, $\Rz$ depends on the transmission and progression closures only through their values at the single point $E_0$, not through their behavior elsewhere and not through their derivatives, and on the uninfected part of the closure only through the location $z^\ast$ of the disease-free equilibrium.
\end{theorem}

\begin{proof}
By \eqref{eq:Frank1}, $K(\phi)=ab^{\!\top}$ with $a=z^\ast_{s_\ast}w\ge0$ and $b^{\!\top}=(\beta^\phi(E_0))^{\!\top}V(\phi)^{-1}\ge0$. A rank-one matrix $ab^{\!\top}$ has eigenvalue $b^{\!\top}a$ with eigenvector $a$ and eigenvalue $0$ on the hyperplane $b^{\perp}$, and $b^{\!\top}a\ge0$, so $\rho(K)=b^{\!\top}a$, which is \eqref{eq:R0closed}. For the last statement, the entries of $V(\phi)$ contributed by closure fluxes $x_j\G_{\phi,k}(x)$ with $j\in\Ical$ equal $\G_{\phi,k}(E_0)$, the derivative terms again being multiplied by $x_j=0$; fluxes of the uninfected part have no infected component and affect \eqref{eq:R0closed} only through $z^\ast$.
\end{proof}

Formula \eqref{eq:R0closed} has the classical interpretation: $V^{-1}w$ is the vector of expected sojourn times in the infected classes for an individual entering the infected pool, and $\beta^\phi(E_0)$ gives the transmission rates per unit time at the prevailing susceptible density $z^\ast_{s_\ast}$. Only the value of the network at $E_0$ enters: $\Rz$ is computable exactly for a trained model by one forward evaluation and one linear solve, and it is insensitive to everything the closure has learned away from $E_0$; Section \ref{sec:identifiability} draws the consequences.

\begin{theorem}[Monotone bounds]\label{thm:R0bounds}
Let $F_1,F_2\ge0$ with $F_1\le F_2$ and let $V$ be a nonsingular M matrix. Then $\rho(F_1V^{-1})\le\rho(F_2V^{-1})$. Consequently, if $\beta^{-}\le\beta^{\phi}(E_0)\le\beta^{+}$ then
\begin{equation}\label{eq:R0sandwich}
z^{\ast}_{s_\ast}(\beta^{-})^{\!\top}V^{-1}w\ \le\ \Rz(\phi)\ \le\ z^{\ast}_{s_\ast}(\beta^{+})^{\!\top}V^{-1}w,
\end{equation}
and if $\phi_1,\phi_2$ produce the same $V$ and the same $z^\ast$ and $\|\beta^{\phi_1}(E_0)-\beta^{\phi_2}(E_0)\|_\infty\le\varepsilon$, then
\begin{equation}\label{eq:R0lip}
|\Rz(\phi_1)-\Rz(\phi_2)|\le \varepsilon\,z^\ast_{s_\ast}\,\one^{\!\top}V^{-1}w .
\end{equation}
\end{theorem}

\begin{proof}
Since $V^{-1}\ge0$ by Proposition \ref{prop:ngm}, $0\le F_1V^{-1}\le F_2V^{-1}$. For $0\le A\le B$ one has $0\le A^k\le B^k$ for every $k$, hence $\|A^k\|_\infty\le\|B^k\|_\infty$ for the maximum absolute row sum norm, and Gelfand's formula gives $\rho(A)=\lim_k\|A^k\|_\infty^{1/k}\le\rho(B)$ \cite{HornJohnson2013}. Inequalities \eqref{eq:R0sandwich} and \eqref{eq:R0lip} then follow from \eqref{eq:R0closed} and $V^{-1}w\ge0$.
\end{proof}

\begin{remark}\label{rem:Vpert}
If the closure also contributes to $V$, then $\Rz$ is no longer linear in the learned quantities but remains smooth: $K=ab^{\!\top}$ has rank at most one, so whenever $\Rz=b^{\!\top}a>0$ it is a simple eigenvalue with right eigenvector $\omega=a$ and left eigenvector $\upsilon=b$, hence analytic in the entries of $K$, and differentiating $K\omega=\rho\omega$ and multiplying by $\upsilon^{\!\top}$, which annihilates the term in $\partial\omega$ because $\upsilon^{\!\top}K=\rho\upsilon^{\!\top}$, gives
\begin{equation}\label{eq:perron-deriv}
\frac{\partial\rho(K)}{\partial K_{ij}}=\frac{\upsilon_i\omega_j}{\upsilon^{\!\top}\omega}.
\end{equation}
\end{remark}

\begin{example}\label{ex:R0seir}
For \eqref{eq:seir}, $y=(E,I)^{\!\top}$, $w=(1,0)^{\!\top}$ and $\beta^\phi(E_0)=(0,\beta^\ast)^{\!\top}$ with $\beta^\ast:=\beta_\phi(E_0)$, so that
\[
V=\begin{pmatrix}\mu+\sigma&0\\-\sigma&\mu+\gamma+\delta\end{pmatrix},\qquad
V^{-1}w=\frac{1}{(\mu+\sigma)(\mu+\gamma+\delta)}\begin{pmatrix}\mu+\gamma+\delta\\ \sigma\end{pmatrix},
\]
and \eqref{eq:R0closed} gives
\begin{equation}\label{eq:R0seir}
\Rz(\phi)=\frac{\sigma\beta^{\ast}\Lambda}{\mu(\mu+\sigma)(\mu+\gamma+\delta)} .
\end{equation}
This is the classical expression with the transmission coefficient replaced by the value of the learned function at $E_0$; no derivative of $\beta_\phi$ appears. Since $V$ is lower triangular with positive diagonal, $s(-V)=-\min\{\mu+\sigma,\mu+\gamma+\delta\}<0$, so Assumption \ref{as:V} holds for every $\phi$.
\end{example}

\section{Endemic equilibria}\label{sec:endemic}

Since $\beta^\phi$ is a black box, no algebraic elimination is available, and we reduce the equilibrium problem to a scalar fixed point equation for the force of infection. The reduction is what makes the endemic theory independent of the particular closure: existence above threshold (Theorem \ref{thm:EEexist}) uses only the slope $H'(0)=\Rz$ and an a priori bound on $H$, which hold for every admissible closure, whereas uniqueness (Theorem \ref{thm:EEunique}) needs the shape condition \ref{cond:sublin}, and the branch hypothesis, Assumption \ref{as:A6}, has to be checked model by model (Remark \ref{rem:A6examples}). Freezing the force of infection at a constant value $\lambda$ turns the transmission term into a linear inflow; the equilibria $x(\lambda)$ of the frozen system, traced as $\lambda$ varies, form the branch along which endemic equilibria are sought, and $\lambda$ is the branch parameter. On $\Om$ one has $\lambda_\phi(x)\le\bar G\sum_{j\in\Ical}x_j\le\bar G\Lambda/\mu=:\lambda^\sharp$; fix $\bar\lambda>\lambda^\sharp$.

\begin{assumption}[Conditioned branch]\label{as:A6}
For $\lambda\in[0,\bar\lambda]$ let $F_\phi^{[\lambda]}$ be obtained from $F_\phi$ by replacing $\lambda_\phi(x)$ in \eqref{eq:foi} by the constant $\lambda$. For each such $\lambda$ the equation $F^{[\lambda]}_\phi(x)=0$ has a unique solution $x(\lambda)\in\Om$, and $\lambda\mapsto x(\lambda)$ is continuous with $x(0)=E_0$.
\end{assumption}

\begin{remark}\label{rem:A6examples}
Assumption \ref{as:A6} holds whenever the closure has no progression or uninfected part and $f(x,\theta)=b+Mx$ is affine. Indeed, Assumption \ref{as:A2} gives $b\ge0$ and $M$ with nonnegative off diagonal entries, Assumption \ref{as:A3} gives $\one^{\!\top}b\le\Lambda$ and $\one^{\!\top}M\le-\mu\one^{\!\top}$, and the frozen transmission term adds $\lambda(\hat w-e_{s_\ast})e_{s_\ast}^{\!\top}$, which has nonnegative off diagonal entries and zero column sums. Hence $-M_\lambda$, with $M_\lambda$ the resulting matrix, is a nonsingular M matrix, the unique equilibrium $x(\lambda)=-M_\lambda^{-1}b\ge0$ depends continuously on $\lambda$, and summing the equilibrium equations gives $\one^{\!\top}x(\lambda)\le\one^{\!\top}b/\mu\le\Lambda/\mu$, so $x(\lambda)\in\Om$. For \eqref{eq:seir} the frozen equilibrium equations give explicitly
\begin{equation}\label{eq:seirbranch}
S(\lambda)=\frac{\Lambda}{\lambda+\mu},\qquad E(\lambda)=\frac{\lambda S(\lambda)}{\mu+\sigma},\qquad I(\lambda)=\frac{\sigma E(\lambda)}{\mu+\gamma+\delta},\qquad R(\lambda)=\frac{\gamma I(\lambda)}{\mu},
\end{equation}
a continuous branch with $x(0)=E_0$ and $\one^{\!\top}x(\lambda)=(\Lambda-\delta I(\lambda))/\mu\le\Lambda/\mu$, so $x(\lambda)\in\Om$; for \eqref{eq:sir} one drops the $E$ equation and puts $I(\lambda)=\lambda S(\lambda)/(\mu+\gamma+\delta)$. In both cases $H(\lambda)=\lambda_\phi(x(\lambda))$ is explicit. If, in \eqref{eq:seir}, $\beta_\phi=\beta_\phi(I)$ depends only on prevalence, then
\[
\frac{H(\lambda)}{\lambda}=\beta_\phi\bigl(I(\lambda)\bigr)\,\frac{\sigma\Lambda}{(\mu+\sigma)(\mu+\gamma+\delta)(\lambda+\mu)},
\]
and since $I(\lambda)$ is increasing in $\lambda$ while $(\lambda+\mu)^{-1}$ is strictly decreasing, Condition \ref{cond:sublin} below holds whenever $\beta_\phi$ is positive and nonincreasing in $I$.
\end{remark}

\begin{lemma}\label{lem:branch}
Under Assumptions \ref{as:A1}--\ref{as:A5} and \ref{as:V} there are $\lambda_0>0$ and a neighborhood $U$ of $E_0$ such that for $\lambda\in[0,\lambda_0)$ the frozen equation $F^{[\lambda]}_\phi(x)=0$ has a unique solution $x(\lambda)\in U$; it is $C^1$ in $\lambda$, $x(0)=E_0$, and its infected component satisfies $y'(0)=z^\ast_{s_\ast}V(\phi)^{-1}w\ge0$. If Assumption \ref{as:A6} holds, its branch coincides with this local branch for small $\lambda$.
\end{lemma}

Lemma \ref{lem:branch} is local: it does not give existence or uniqueness of frozen equilibria in $\Om$ for larger $\lambda$, which is why Assumption \ref{as:A6} is imposed separately.

\begin{proof}
At $(\lambda,x)=(0,E_0)$ the equilibrium equations hold. Their Jacobian in $x$ is block triangular: the derivative of the infected equations with respect to $z$ vanishes because $\Fcal^{[0]}_\phi$ and $\Vcal_\phi$ vanish identically on $\mathcal X_0$ (Proposition \ref{prop:dfinv}); the infected diagonal block is $-V(\phi)$, since $D_y\Fcal^{[\lambda]}_\phi=0$ for frozen $\lambda$; and the uninfected diagonal block is $Dg_\phi(z^\ast)$. Both are invertible by Proposition \ref{prop:ngm} and Assumption \ref{as:A5}, so the implicit function theorem applies. Differentiating the infected equations in $\lambda$ at $0$ gives $wz^\ast_{s_\ast}-V(\phi)y'(0)=0$. If Assumption \ref{as:A6} holds, its branch is continuous with value $E_0$ at $\lambda=0$, so it lies in $U$ for small $\lambda$ and coincides there with the solution given by the implicit function theorem.
\end{proof}

\begin{lemma}\label{lem:Hprops}
Under Assumptions \ref{as:A1}--\ref{as:A6} and \ref{as:V}, define $H(\lambda):=\lambda_\phi(x(\lambda))$ on $[0,\bar\lambda]$. Then $H$ is continuous, $H(0)=0$, $0\le H\le\lambda^\sharp<\bar\lambda$, and $H$ is differentiable at $0$ with $H'(0)=\Rz(\phi)$. A point $\bar x\in\Om\setminus\{E_0\}$ is an equilibrium of \eqref{eq:ude} if and only if $\bar x=x(\lambda)$ for some $\lambda\in(0,\bar\lambda]$ with $H(\lambda)=\lambda$.
\end{lemma}

\begin{proof}
Continuity and the bounds are immediate from Assumption \ref{as:A6} and $\lambda_\phi\le\lambda^\sharp$ on $\Om$; $H(0)=0$ because $x(0)=E_0$ has $y=0$. Differentiating $H(\lambda)=(\beta^\phi(x(\lambda)))^{\!\top}y(\lambda)$ at $\lambda=0$ eliminates the terms with derivatives of $\beta^\phi$ because $y(0)=0$, leaving $H'(0)=(\beta^\phi(E_0))^{\!\top}y'(0)=\Rz(\phi)$ by Lemma \ref{lem:branch} and \eqref{eq:R0closed}. If $\bar x$ is an equilibrium of \eqref{eq:ude}, it solves $F^{[\lambda]}_\phi(\bar x)=0$ with $\lambda:=\lambda_\phi(\bar x)$, hence $\bar x=x(\lambda)$ and $H(\lambda)=\lambda$; conversely $H(\lambda)=\lambda$ gives $F_\phi(x(\lambda))=0$. Finally $\lambda=0$ corresponds to $E_0$: an equilibrium $\bar x$ with $\lambda_\phi(\bar x)=0$ solves $F^{[0]}_\phi(\bar x)=0$, hence $\bar x=x(0)=E_0$ by the uniqueness in Assumption \ref{as:A6}.
\end{proof}

\begin{theorem}[Existence]\label{thm:EEexist}
If $\Rz(\phi)>1$, then under Assumptions \ref{as:A1}--\ref{as:A6} and \ref{as:V} the system \eqref{eq:ude} has at least one endemic equilibrium $E_\ast=x(\lambda_\ast)\in\Om$ with $\lambda_\ast\in(0,\lambda^\sharp]$ and $y(\lambda_\ast)\gneq0$.
\end{theorem}

\begin{proof}
Let $\varkappa(\lambda):=H(\lambda)-\lambda$, continuous on $[0,\bar\lambda]$ with $\varkappa(0)=0$ and $\varkappa'(0)=\Rz(\phi)-1>0$, so $\varkappa(\lambda_1)>0$ for some $\lambda_1\in(0,\bar\lambda)$. Since $H\le\lambda^\sharp<\bar\lambda$ by Lemma \ref{lem:Hprops}, $\varkappa(\bar\lambda)<0$. By the intermediate value theorem there is $\lambda_\ast\in(\lambda_1,\bar\lambda)$ with $H(\lambda_\ast)=\lambda_\ast$, and then $\lambda_\ast=H(\lambda_\ast)\le\lambda^\sharp$. Lemma \ref{lem:Hprops} identifies $E_\ast:=x(\lambda_\ast)$ as an equilibrium, and $\lambda_\phi(E_\ast)=\lambda_\ast>0$ forces $y(\lambda_\ast)\gneq0$.
\end{proof}

\begin{condition}[Sublinearity]\label{cond:sublin}
The map $\lambda\mapsto H(\lambda)/\lambda$ is strictly decreasing on $(0,\bar\lambda]$.
\end{condition}

\begin{theorem}[Uniqueness and threshold dichotomy]\label{thm:EEunique}
Under Assumptions \ref{as:A1}--\ref{as:A6}, \ref{as:V} and Condition \ref{cond:sublin}, \eqref{eq:ude} has exactly one endemic equilibrium in $\Om$ if $\Rz(\phi)>1$ and none if $\Rz(\phi)\le1$.
\end{theorem}

\begin{proof}
Endemic equilibria correspond to solutions of $H(\lambda)/\lambda=1$ in $(0,\bar\lambda]$, and the left hand side is strictly decreasing, so there is at most one. By Lemma \ref{lem:Hprops}, $H(\lambda)/\lambda\to\Rz(\phi)$ as $\lambda\downarrow0$, so strict monotonicity gives $H(\lambda)/\lambda<\Rz(\phi)$ for $\lambda>0$; if $\Rz(\phi)\le1$ there is no solution. If $\Rz(\phi)>1$, Theorem \ref{thm:EEexist} gives existence.
\end{proof}

Condition \ref{cond:sublin} is the abstract form of concavity of the incidence in the infection pressure \cite{CapassoSerio1978,Korobeinikov2007}; it fails in particular in the regimes producing multiple endemic states, analyzed in Section \ref{sec:bifurcation}. In $\Ecal_{\mathrm{sep}}$ it reduces to a transparent architectural requirement.

\begin{proposition}\label{prop:EEsep}
Consider \eqref{eq:sir} with $h_\phi\in C^1$, $h_\phi(0)=0$, $h_\phi(I)>0$ for $I>0$, $h_\phi$ nondecreasing, and $I\mapsto h_\phi(I)/I$ nonincreasing on $(0,\Lambda/\mu]$, and set $\Rz(\phi)=\Lambda h'_\phi(0)/(\mu(\mu+\gamma+\delta))$. If $\Rz(\phi)>1$ there is exactly one endemic equilibrium $E_\ast=(S^\ast,I^\ast,R^\ast)$ in $\Om$, with all components positive and
\begin{equation}\label{eq:EEsep}
S^{\ast}=\frac{(\mu+\gamma+\delta)I^{\ast}}{h_\phi(I^{\ast})},\qquad R^{\ast}=\frac{\gamma I^{\ast}}{\mu},
\end{equation}
where $I^\ast$ is the unique zero in $(0,\Lambda/(\mu+\gamma+\delta))$ of
\begin{equation}\label{eq:Theta}
\Pi(I):=\frac{\mu(\mu+\gamma+\delta)I}{h_\phi(I)}+(\mu+\gamma+\delta)I-\Lambda .
\end{equation}
If $\Rz(\phi)\le1$ no endemic equilibrium exists.
\end{proposition}

\begin{proof}
At an equilibrium with $I>0$ the second and third equations of \eqref{eq:sir} give \eqref{eq:EEsep}, and substitution into $\Lambda-S^\ast h_\phi(I^\ast)-\mu S^\ast=0$ gives $\Pi(I^\ast)=0$. Since $h_\phi(I)/I$ is nonincreasing, $I/h_\phi(I)$ is nondecreasing, so $\Pi$ is continuous and strictly increasing on $(0,I^\dagger]$ with $I^\dagger:=\Lambda/(\mu+\gamma+\delta)$. As $h_\phi(I)/I\to h'_\phi(0)$ when $I\downarrow0$, we get $\Pi(0^{+})=\mu(\mu+\gamma+\delta)/h'_\phi(0)-\Lambda$, which is negative exactly when $\Rz(\phi)>1$, while $\Pi(I^{\dagger})=\mu(\mu+\gamma+\delta)I^\dagger/h_\phi(I^\dagger)>0$ and, more generally, $\Pi(I)>0$ for $I\ge I^\dagger$, since then $(\mu+\gamma+\delta)I\ge\Lambda$, so every zero of $\Pi$ lies in $(0,I^\dagger)$. Strict monotonicity gives a unique zero when $\Rz(\phi)>1$ and none otherwise. Membership in $\Om$ follows from Theorem \ref{thm:region}.
\end{proof}

\section{Stability}\label{sec:stability}

\begin{lemma}[M matrix threshold lemma]\label{lem:Mmatrix}
Let $F\ge0$ and let $V$ be a nonsingular M matrix. Then $s(F-V)<0$ if and only if $\rho(FV^{-1})<1$, $s(F-V)=0$ if and only if $\rho(FV^{-1})=1$, and $s(F-V)>0$ if and only if $\rho(FV^{-1})>1$.
\end{lemma}

\begin{proof}
Let $A:=F-V$, which has nonnegative off diagonal entries, choose $s_0>\max_i|A_{ii}|$ and put $B:=A+s_0I\ge0$. Since $\rho(B)\in\spec(B)$ by the Perron-Frobenius theorem and $\spec(A)=\spec(B)-s_0$, we have $s(A)=\rho(B)-s_0$.

Assume $\rho(FV^{-1})<1$. As $FV^{-1}$ and $V^{-1}F$ have the same nonzero eigenvalues, $\rho(V^{-1}F)<1$, so $(I-V^{-1}F)^{-1}=\sum_{k\ge0}(V^{-1}F)^k\ge0$ and $(V-F)^{-1}=(I-V^{-1}F)^{-1}V^{-1}\ge0$. If $\rho(B)=s_0$ then $V-F=s_0I-B$ would be singular; if $\rho(B)>s_0$, take $\xi\gneq0$ with $B\xi=\rho(B)\xi$, so $\xi=(s_0-\rho(B))(V-F)^{-1}\xi\le0$, forcing $\xi=0$. Hence $\rho(B)<s_0$ and $s(A)<0$.

Assume conversely $\rho:=\rho(FV^{-1})\ge1$ and take $u\gneq0$ with $FV^{-1}u=\rho u$. Put $p:=V^{-1}u\ge0$, $p\ne0$. Then $Fp=\rho Vp$, so $Ap=(\rho-1)Vp=(\rho-1)u\ge0$ and $Bp\ge s_0p$. Since $B\ge0$ preserves the order, $B^kp\ge s_0^kp$ for all $k$, so $\|B^k\|\ge s_0^k$ and $\rho(B)=\lim_k\|B^k\|^{1/k}\ge s_0$, that is $s(A)\ge0$.

It remains to separate the strict cases. If $\rho(FV^{-1})\le1$, then $\rho(tFV^{-1})<1$ for $t\in(0,1)$, so $s(tF-V)<0$ by the first part, and continuity of the spectrum gives $s(F-V)\le0$. If $\rho(FV^{-1})>1$ but $s(F-V)\le0$, then for every $\epsilon>0$ the matrix $V+\epsilon I$ is a nonsingular M matrix with $s(F-V-\epsilon I)<0$, so $\rho(F(V+\epsilon I)^{-1})<1$ by the first part, and letting $\epsilon\downarrow0$ gives $\rho(FV^{-1})\le1$, a contradiction. Hence $s(F-V)>0$ if and only if $\rho(FV^{-1})>1$, and the case of equality follows.
\end{proof}

\subsection{Local stability}

\begin{theorem}\label{thm:LASDFE}
Under Assumptions \ref{as:A1}--\ref{as:A5} and \ref{as:V} the Jacobian of \eqref{eq:ude} at $E_0$ is block triangular, with the infected coordinates listed first,
\begin{equation}\label{eq:J0}
J(E_0)=\begin{pmatrix}F(\phi)-V(\phi)&0\\ \ast&Dg_\phi(z^\ast)\end{pmatrix},
\end{equation}
so $\spec(J(E_0))=\spec(F(\phi)-V(\phi))\cup\spec(Dg_\phi(z^\ast))$. Hence $E_0$ is locally asymptotically stable if $\Rz(\phi)<1$ and unstable if $\Rz(\phi)>1$.
\end{theorem}

\begin{proof}
Write the infected equations as $\dot y=\Psi(x)$ with $\Psi:=\Fcal_\phi-\Vcal_\phi$. By Proposition \ref{prop:dfinv}, $\Psi(z,0)=0$ for all $z$; differentiating this identity in $z$ gives $D_z\Psi(E_0)=0$, while $D_y\Psi(E_0)=F(\phi)-V(\phi)$ by \eqref{eq:FVmat}. The uninfected block is $Dg_\phi(z^\ast)$, which has $s<0$ by Assumption \ref{as:A5}, so the stability type is decided by $s(F-V)$, whose sign is given by Lemma \ref{lem:Mmatrix}. The conclusion follows from the principle of linearized stability \cite[Ch.~IX]{Hartman2002}.
\end{proof}

\begin{theorem}\label{thm:LASEE}
Under the hypotheses of Proposition \ref{prop:EEsep} with $\Rz(\phi)>1$, the endemic equilibrium $E_\ast$ of \eqref{eq:sir} is locally asymptotically stable.
\end{theorem}

\begin{proof}
The $R$ equation is linear and decoupled, contributing the eigenvalue $-\mu$, so it suffices to consider the $(S,I)$ system, whose Jacobian at $E_\ast$ is
\[
J_\ast=\begin{pmatrix}-h_\phi(I^\ast)-\mu&-S^\ast h'_\phi(I^\ast)\\[2pt] h_\phi(I^\ast)&S^\ast h'_\phi(I^\ast)-(\mu+\gamma+\delta)\end{pmatrix}.
\]
Sublinearity gives $(h_\phi(I)/I)'=(Ih'_\phi(I)-h_\phi(I))/I^2\le0$, hence $I^\ast h'_\phi(I^\ast)\le h_\phi(I^\ast)$, and using $S^\ast h_\phi(I^\ast)=(\mu+\gamma+\delta)I^\ast$ we obtain $S^\ast h'_\phi(I^\ast)\le(\mu+\gamma+\delta)$. Therefore $\tr J_\ast\le-h_\phi(I^\ast)-\mu<0$ and
\begin{align*}
\det J_\ast&=\bigl(h_\phi(I^\ast)+\mu\bigr)\bigl(\mu+\gamma+\delta-S^\ast h'_\phi(I^\ast)\bigr)+S^\ast h'_\phi(I^\ast)h_\phi(I^\ast)\\
&=\mu\bigl(\mu+\gamma+\delta-S^\ast h'_\phi(I^\ast)\bigr)+(\mu+\gamma+\delta)h_\phi(I^\ast)>0 .
\end{align*}
The Routh-Hurwitz conditions for a $2\times2$ matrix are exactly $\tr J_\ast<0$ and $\det J_\ast>0$.
\end{proof}

\subsection{Global stability of the disease-free equilibrium}

\begin{condition}[Uniform domination]\label{cond:dom}
For all $x\in\Om$,
\begin{equation}\label{eq:dom}
x_{s_\ast}\beta^\phi(x)\ \le\ z^\ast_{s_\ast}\beta^\phi(E_0),\qquad \Vcal_\phi(x)\ \ge\ V(\phi)\,y .
\end{equation}
\end{condition}

Suppose the susceptible class carries the whole population at $E_0$, that is $E_0=z^\ast_{s_\ast}e_{s_\ast}$ with $z^\ast_{s_\ast}=\Lambda/\mu$, as in \eqref{eq:seir}. Then the first inequality of \eqref{eq:dom} holds whenever $\beta^\phi$ is nonincreasing in every coordinate other than $x_{s_\ast}$ and the map $\vartheta\mapsto\vartheta\,\beta^\phi(\vartheta e_{s_\ast})$ is componentwise nondecreasing on $[0,z^\ast_{s_\ast}]$: for $x\in\Om$ one has $x_{s_\ast}\le z^\ast_{s_\ast}$ and $x\ge x_{s_\ast}e_{s_\ast}$, so $x_{s_\ast}\beta^\phi(x)\le x_{s_\ast}\beta^\phi(x_{s_\ast}e_{s_\ast})\le z^\ast_{s_\ast}\beta^\phi(E_0)$. The requirement on $\vartheta\mapsto\vartheta\,\beta^\phi(\vartheta e_{s_\ast})$ does not force $\beta^\phi$ to be monotone in the susceptible coordinate: it holds when $\beta^\phi$ is nondecreasing in $x_{s_\ast}$, and it also allows saturation in $S$ provided the elasticity $\partial\ln\beta^\phi_j/\partial\ln x_{s_\ast}$ along the susceptible axis does not fall below $-1$, a threshold that reappears in Section \ref{sec:sensitivity}. A dependence on $x_{s_\ast}$ that is nonincreasing, but with elasticity below $-1$ on part of $[0,z^\ast_{s_\ast}]$, can violate \eqref{eq:dom}: for $\beta^\phi(x)=\beta_0\ee^{-\alpha x_{s_\ast}}$ with $\alpha z^\ast_{s_\ast}=3$ the maximum of $x_{s_\ast}\beta^\phi$ on $\Om$ exceeds $z^\ast_{s_\ast}\beta^\phi(E_0)$ by a factor $\ee^{2}/3$. Monotone architectures enforce the sign conditions on the coordinates other than $x_{s_\ast}$ by construction. The second inequality of \eqref{eq:dom} holds with equality when the transfers among infected classes are mechanistic and linear. Both are conditions on the compact set $\Om$ and can be certified after training.

The following lemma identifies the dynamics on the disease-free set; it is used for global stability and again for persistence.

\begin{lemma}\label{lem:boundary}
Under Assumptions \ref{as:A1}, \ref{as:A4}, \ref{as:A5} and Definition \ref{def:closure}, $\{E_0\}$ is the only nonempty compact invariant subset of $\Om\cap\mathcal X_0$.
\end{lemma}

\begin{proof}
Let $M\subset\Om\cap\mathcal X_0$ be nonempty, compact and invariant, and let $p\in M$. Invariance provides, for every $t>0$, a point $p_{-t}\in M$ with $p=x(t;p_{-t})$, so the entire orbit through $p$ lies in $M$ and, by Proposition \ref{prop:dfinv}, is an orbit of the reduced system \eqref{eq:dfsub}. Its alpha limit set $\mathcal A(p)\subset M$ is nonempty, compact and invariant; for $q\in\mathcal A(p)$, invariance and closedness give $\omega(q)\subset\mathcal A(p)$, and $\omega(q)=\{E_0\}$ by the global attractivity in Assumption \ref{as:A5}, so $E_0\in\mathcal A(p)$. Given $\varepsilon>0$, let $\delta>0$ be provided by the Lyapunov stability of $E_0$ for the reduced system; since $E_0\in\mathcal A(p)$ there is $t>0$ with $p_{-t}\in B_\delta(E_0)$, hence $p=x(t;p_{-t})\in B_\varepsilon(E_0)$. As $\varepsilon$ was arbitrary, $p=E_0$.
\end{proof}

\begin{theorem}\label{thm:GASDFE}
Let Assumptions \ref{as:A1}--\ref{as:A5}, \ref{as:V} and Condition \ref{cond:dom} hold. If $\Rz(\phi)<1$ then $E_0$ is globally asymptotically stable in $\Om$, and for every $\alpha\in(0,-s(F(\phi)-V(\phi)))$ there is $C_\alpha>0$, depending only on $F(\phi)-V(\phi)$ and $\alpha$, with
\begin{equation}\label{eq:expdecay}
\|y(t)\|\le C_\alpha\ee^{-\alpha t}\|y(0)\|,\qquad t\ge0 ;
\end{equation}
here $-s(F(\phi)-V(\phi))>0$ by Lemma \ref{lem:Mmatrix}. If $\Rz(\phi)>1$ then $E_0$ is unstable.
\end{theorem}

\begin{proof}
Let $x_0\in\Om$; the solution stays in $\Om$ by Theorem \ref{thm:region}. Because $w\ge0$ and $y\ge0$, Condition \ref{cond:dom} gives $\Fcal_\phi(x)=w\,x_{s_\ast}(\beta^\phi(x))^{\!\top}y\le w\,z^\ast_{s_\ast}(\beta^\phi(E_0))^{\!\top}y=Fy$ and $\Vcal_\phi(x)\ge Vy$, so the infected block satisfies the componentwise differential inequality $\dot y\le Ay$ with $A:=F-V$. The matrix $A$ has nonnegative off diagonal entries, so the linear system $\dot\eta=A\eta$ is cooperative and the Kamke comparison theorem \cite{Smith1995,Walter1998} yields $0\le y(t)\le\ee^{At}y(0)$ for $t\ge0$. Since $\rho(FV^{-1})=\Rz<1$, Lemma \ref{lem:Mmatrix} gives $s(A)<0$, and for every $\alpha\in(0,-s(A))$ there is $C_\alpha$ with $\|\ee^{At}\|\le C_\alpha\ee^{-\alpha t}$; this is \eqref{eq:expdecay}, using $\|y(t)\|\le\|\ee^{At}y(0)\|$ for $0\le y(t)\le\ee^{At}y(0)$.

In particular $y(t)\to0$, so the omega limit set $\omega(x_0)$, which is nonempty, compact and invariant because the solution is bounded, lies in $\Om\cap\mathcal X_0$, and Lemma \ref{lem:boundary} gives $\omega(x_0)=\{E_0\}$, that is $x(t)\to E_0$ for every $x_0\in\Om$. Together with the local asymptotic stability in Theorem \ref{thm:LASDFE} this is global asymptotic stability in $\Om$; instability for $\Rz>1$ is Theorem \ref{thm:LASDFE}.
\end{proof}

\begin{remark}\label{rem:matrixlyap}
Theorem \ref{thm:GASDFE} is the matrix-theoretic method of Shuai and van den Driessche \cite{ShuaiVdD2013} adapted to a learned closure: there, under the same domination condition, $L(x):=u^{\!\top}V^{-1}y$ with $u$ a left Perron vector of $V^{-1}F$ is a Lyapunov function, $\dot L\le(\Rz-1)u^{\!\top}y\le0$. When $V^{-1}F$ is irreducible one has $u>0$, hence $\varpi^{\!\top}:=u^{\!\top}V^{-1}>0$ by Proposition \ref{prop:ngm}, and the argument gives the explicit rate $\|y(t)\|\le C\ee^{-c(1-\Rz)t}\|y(0)\|$ with $c=\min_iu_i/\varpi_i$ and $C=\sqrt m\max_i\varpi_i/\min_i\varpi_i$. Irreducibility fails, however, whenever some infected class is not infectious, that is $\beta^\phi_j(E_0)=0$ for some $j$: for \eqref{eq:seir} one has $\beta^\phi(E_0)=(0,\beta^\ast)^{\!\top}$, so $V^{-1}F$ has a zero column and $u=(0,\beta^\ast)^{\!\top}$, and the constant $c$ vanishes. The comparison argument above needs no irreducibility and covers this case; the vector $\varpi^{\!\top}=(\beta^\phi(E_0))^{\!\top}V^{-1}$ is nevertheless strictly positive for \eqref{eq:seir}, so $L$ remains a Lyapunov function there.
\end{remark}

\subsection{Global stability of the endemic equilibrium}

Without constraints on the closure the endemic state need not be globally stable, since models with nonmonotone incidence admit periodic and bistable behavior \cite{LiuLevinIwasa1986,RuanWang2003}. In $\Ecal_{\mathrm{sep}}$ the classical conditions transfer verbatim to the network.

\begin{theorem}\label{thm:GASEE}
Consider \eqref{eq:sir} under the hypotheses of Proposition \ref{prop:EEsep} with $\Rz(\phi)>1$. Then the endemic equilibrium $E_\ast$ is globally asymptotically stable in $\{x\in\Om:I>0\}$.
\end{theorem}

\begin{proof}
Since $\dot R=\gamma I-\mu R$ is decoupled and linear, $I(t)\to I^\ast$ implies $R(t)\to R^\ast$, so it suffices to treat the $(S,I)$ system on $\Om_2:=\{(S,I)\in\Rp^2:S+I\le\Lambda/\mu\}$, which is compact and positively invariant by Theorem \ref{thm:pos} and $\frac{\dd}{\dd t}(S+I)\le\Lambda-\mu(S+I)$; moreover $I(t)>0$ for $t\ge0$, since $\dot I\ge-(\mu+\gamma+\delta)I$, and $S(t)>0$ for $t>0$, since $\dot S=\Lambda>0$ whenever $S=0$. Define, on the interior of $\Om_2$,
\begin{equation}\label{eq:Vee}
W(S,I):=S-S^{\ast}-S^{\ast}\ln\frac{S}{S^{\ast}}+I-I^{\ast}-I^{\ast}\ln\frac{I}{I^{\ast}},
\end{equation}
which is $C^1$, positive definite with respect to $E_\ast$, and satisfies $W(S,I)\to\infty$ as $S\downarrow0$ or $I\downarrow0$. Consequently every sublevel set $\{W\le W(x_0)\}\cap\Om_2$ is a compact subset of the interior of $\Om_2$, so LaSalle's principle applies there once $\dot W\le0$ is established.

Using $\Lambda=\mu S^\ast+S^\ast h_\phi(I^\ast)$ and $\mu+\gamma+\delta=S^\ast h_\phi(I^\ast)/I^\ast$,
\[
\dot W=\Bigl(1-\frac{S^\ast}{S}\Bigr)\bigl(\Lambda-Sh_\phi(I)-\mu S\bigr)+\Bigl(1-\frac{I^\ast}{I}\Bigr)\bigl(Sh_\phi(I)-(\mu+\gamma+\delta)I\bigr)
=-\frac{\mu(S-S^\ast)^2}{S}+S^\ast h_\phi(I^\ast)\,\mathcal E,
\]
where, writing $s:=S/S^\ast$, $\iota:=I/I^\ast$ and $\zeta:=h_\phi(I)/h_\phi(I^\ast)$, expansion of the remaining terms gives
\[
S^\ast h_\phi(I^\ast)\mathcal E=2S^\ast h_\phi(I^\ast)-\frac{S^{\ast2}h_\phi(I^\ast)}{S}+S^\ast h_\phi(I)-\frac{S^\ast h_\phi(I^\ast)I}{I^\ast}-\frac{I^\ast Sh_\phi(I)}{I},
\qquad \mathcal E=2-\frac1s+\zeta-\iota-\frac{s\zeta}{\iota}.
\]
Since $\zeta-\iota-1+\iota/\zeta=(\zeta-1)(\zeta-\iota)/\zeta$, we have the identity
\begin{equation}\label{eq:decomp}
\mathcal E=\Bigl(3-\frac1s-\frac{s\zeta}{\iota}-\frac{\iota}{\zeta}\Bigr)+\frac{(\zeta-1)(\zeta-\iota)}{\zeta}.
\end{equation}
The three positive numbers $1/s$, $s\zeta/\iota$ and $\iota/\zeta$ have product one, so their sum is at least three by the arithmetic geometric mean inequality, and the first bracket is nonpositive, vanishing only if $s=1$ and $\zeta=\iota$. For the second term, monotonicity of $h_\phi$ gives $\zeta\ge1$ for $I\ge I^\ast$ and $\zeta\le1$ for $I\le I^\ast$, while $h_\phi(I)/I$ nonincreasing gives $\zeta\le\iota$ for $I\ge I^\ast$ and $\zeta\ge\iota$ for $I\le I^\ast$; in both cases the factors have opposite signs and the term is nonpositive. Hence $\dot W\le0$, and $\dot W=0$ requires $S=S^\ast$ and $\zeta=\iota$. On any invariant subset of $\{\dot W=0\}$ we have $S\equiv S^\ast$, so $\dot S\equiv0$ gives $S^\ast h_\phi(I)=\Lambda-\mu S^\ast=S^\ast h_\phi(I^\ast)$, that is $\zeta=1$, and with $\zeta=\iota$ this yields $I=I^\ast$. The largest invariant subset of $\{\dot W=0\}$ is therefore $\{E_\ast\}$, and LaSalle's principle \cite{LaSalle1976} gives convergence of every solution with $I(0)>0$; combined with Theorem \ref{thm:LASEE} this proves global asymptotic stability.
\end{proof}

\begin{remark}\label{rem:archglobal}
The hypotheses of Theorem \ref{thm:GASEE}, monotonicity of $h_\phi$ and sublinearity of $h_\phi(I)/I$, are the classical conditions for global stability with general incidence \cite{KorobeinikovMaini2004,Korobeinikov2007}; here they are constraints on a network and can be imposed architecturally, monotonicity by nonnegative weights with monotone activations and sublinearity by the representation $h_\phi(I)=I\varrho(\N_\phi(I))$ with $\varrho$ nonnegative and nonincreasing in a monotone architecture. The latter representation does not by itself make $h_\phi$ monotone; imposing both properties at once requires in addition that the elasticity of $I\mapsto\varrho(\N_\phi(I))$ be at least $-1$, a condition that can be imposed as a constraint or certified on $[0,\Lambda/\mu]$. A trained model in this class inherits the entire classical global theory.
\end{remark}

\section{Persistence and extinction}\label{sec:persistence}

We use the persistence theorem of Hale and Waltman \cite{HaleWaltman1989}, in the following form. Let $X$ be a complete metric space that is the closure of an open subset $X^{\circ}$, with nonempty boundary $\partial X^{\circ}=X\setminus X^{\circ}$, and let $T(t)$ be a continuous semiflow on $X$ with $T(t)X^{\circ}\subset X^{\circ}$ and $T(t)\partial X^{\circ}\subset\partial X^{\circ}$ for $t\ge0$. Suppose that $T(t)$ is asymptotically smooth and point dissipative, that positive orbits of bounded sets are bounded, and that the maximal compact invariant set $A_\partial$ of $T(t)$ in $\partial X^{\circ}$ is isolated and has an acyclic covering $\{M_1,\dots,M_k\}$ by pairwise disjoint compact isolated invariant sets. Then $T(t)$ is uniformly persistent, meaning that there is $\varepsilon_0>0$ with $\liminf_{t\to\infty}\dist(T(t)x,\partial X^{\circ})\ge\varepsilon_0$ for every $x\in X^{\circ}$, if and only if $W^s(M_i)\cap X^{\circ}=\emptyset$ for each $i$, where $W^s(M_i)$ is the set of points whose omega limit set is contained in $M_i$. The same conclusion is available from the formulations in \cite{Thieme1993,SmithThieme2011,Zhao2017}. In the present setting the spaces are
\begin{equation}\label{eq:persistsets}
X:=\Om,\qquad X^{\circ}:=\{x\in\Om:\ \one^{\!\top}y>0\},\qquad \partial X^{\circ}=\Om\cap\mathcal X_0 ,
\end{equation}
so that the extinction boundary $\partial X^{\circ}$ is the part of the disease-free set \eqref{eq:dfset} inside the feasible region.

\begin{theorem}[Uniform persistence]\label{thm:persist}
Let Assumptions \ref{as:A1}--\ref{as:A5} and \ref{as:V} hold, let $F(\phi)-V(\phi)$ be irreducible and let $\Rz(\phi)>1$. Then there is $\varepsilon_0>0$, independent of the initial datum, such that every solution with $x_0\in\Om$ and $\one^{\!\top}y(0)>0$ satisfies
\begin{equation}\label{eq:persist}
\liminf_{t\to\infty}\one^{\!\top}y(t)\ \ge\ \varepsilon_0 .
\end{equation}
Since solutions are also ultimately bounded by Theorem \ref{thm:region}, the infected population is permanent: it is eventually confined to $[\varepsilon_0/2,\Lambda/\mu]$ in $\ell^1$ norm. If in addition Assumption \ref{as:A6} holds, an endemic equilibrium exists by Theorem \ref{thm:EEexist}.
\end{theorem}

\begin{proof}
We verify the hypotheses of the persistence theorem for the sets \eqref{eq:persistsets} and the semiflow $T(t)x_0:=x(t;x_0)$.

\emph{Semiflow and topology.} By Theorems \ref{thm:local} and \ref{thm:region} the solution map is a continuous semiflow on the compact set $\Om$. The set $X^{\circ}$ is open in $\Om$, its closure is $\Om$ because every point of $\Om\cap\mathcal X_0$ is a limit of points with $\one^{\!\top}y>0$ in $\Om$, and its boundary in $\Om$ is $\Om\cap\mathcal X_0\neq\emptyset$.

\emph{Invariance of $X^{\circ}$ and $\partial X^{\circ}$.} Put $Y:=\one^{\!\top}y$. The argument in the last paragraph of the proof of Theorem \ref{thm:pos}, applied on the box $\{x\in\Rp^n:x\le(\Lambda/\mu)\one\}\supset\Om$ with Lipschitz constant $L'$, gives $F_{\phi,j}(x)\ge-L'x_j$ for every $j\in\Ical$ and $x\in\Om$; hence $\dot Y\ge-L'Y$ and $Y(t)\ge Y(0)\ee^{-L't}>0$, so $X^{\circ}$ is forward invariant. Forward invariance of $\partial X^{\circ}$ is Proposition \ref{prop:dfinv} together with that of $\Om$.

\emph{Compactness and dissipativity.} Because $\Om$ is compact, $T(t)$ is point dissipative and asymptotically smooth, and positive orbits of bounded sets are bounded; the compact global attractor exists by Theorem \ref{thm:region}.

\emph{Boundary invariant set, isolation and acyclicity.} By Lemma \ref{lem:boundary}, the maximal compact invariant set in $\partial X^{\circ}$ is $A_\partial=\{E_0\}$, and we take the covering $\{M_1\}=\{\{E_0\}\}$. By Theorem \ref{thm:LASDFE} and Lemma \ref{lem:Mmatrix}, the spectrum of the Jacobian at $E_0$ is $\spec(F-V)\cup\spec(Dg_\phi(z^\ast))$, and $\Rz>1$ gives $s(F-V)>0$ while Assumption \ref{as:A5} gives $s(Dg_\phi(z^\ast))<0$; here $E_0$ need not be hyperbolic, since $F-V$ may have eigenvalues on the imaginary axis, but it is an isolated invariant set. Indeed, let $\upsilon>0$ and $\varrho>0$ be as in the next step, so that $\frac{\dd}{\dd t}(\upsilon^{\!\top}y)\ge\tfrac12s(F-V)\upsilon^{\!\top}y$ on $B_\varrho(E_0)\cap\Om$. Along a full orbit contained in $B_\varrho(E_0)\cap\Om$ the quantity $\upsilon^{\!\top}y$ is bounded, which this inequality allows only if $y\equiv0$. The orbit is then a full orbit of the reduced system \eqref{eq:dfsub} near $z^\ast$, which is exponentially stable because $s(Dg_\phi(z^\ast))<0$; along any other full orbit near $z^\ast$, a quadratic Lyapunov function that decays exponentially along solutions grows without bound backward in time. Hence, after decreasing $\varrho$ if necessary, which preserves the growth inequality, $\{E_0\}$ is the maximal invariant set in $B_\varrho(E_0)\cap\Om$. A covering by a single set is acyclic unless $E_0$ is chained to itself by a full orbit in $\partial X^{\circ}$ that is homoclinic to $E_0$; such an orbit together with $E_0$ would be a compact invariant subset of $\partial X^{\circ}$ different from $\{E_0\}$, contradicting Lemma \ref{lem:boundary}.

\emph{No interior stable set.} We show $W^s(\{E_0\})\cap X^{\circ}=\emptyset$. Suppose some solution with $x(0)\in X^{\circ}$ satisfies $x(t)\to E_0$. Write $\Psi:=\Fcal_\phi-\Vcal_\phi$, so $\Psi(z,0)\equiv0$ by Proposition \ref{prop:dfinv} and $\Psi(x)=\int_0^1D_y\Psi(z,\tau y)y\,\dd\tau$. By continuity of $D_y\Psi$ and $D_y\Psi(E_0)=F-V$, for every $\varepsilon>0$ there is $\varrho_\varepsilon>0$ such that $\|x-E_0\|\le\varrho_\varepsilon$ implies $\|\Psi(x)-(F-V)y\|_\infty\le\varepsilon\|y\|_1$. Since $F-V$ has nonnegative off diagonal entries and is irreducible, $F-V+s_0I$ is nonnegative and irreducible for large $s_0$, so the Perron-Frobenius theorem provides $\upsilon>0$ with $\upsilon^{\!\top}(F-V)=s(F-V)\upsilon^{\!\top}$. Whenever $x\in\Om$ and $\|x-E_0\|\le\varrho_\varepsilon$, in particular for $t$ large,
\[
\frac{\dd}{\dd t}(\upsilon^{\!\top}y)=\upsilon^{\!\top}\Psi(x)\ \ge\ s(F-V)\upsilon^{\!\top}y-\varepsilon(\upsilon^{\!\top}\one)\|y\|_1\ \ge\ \Bigl(s(F-V)-\frac{\varepsilon\,\upsilon^{\!\top}\one}{\min_i\upsilon_i}\Bigr)\upsilon^{\!\top}y,
\]
using $\one^{\!\top}y\le(\min_i\upsilon_i)^{-1}\upsilon^{\!\top}y$ for $y\ge0$. Choosing $\varepsilon$ so that the bracket exceeds $s(F-V)/2>0$, and putting $\varrho:=\varrho_\varepsilon$ (the radius used in the isolation step), yields exponential growth of $\upsilon^{\!\top}y$, which is positive at all times by the invariance of $X^{\circ}$, contradicting $y(t)\to0$. Irreducibility of $F-V$ is used only to obtain this strictly positive left eigenvector, here and in the isolation step; for \eqref{eq:seir} it holds because $\sigma>0$ and $\beta^\ast>0$.

\emph{Conclusion.} The persistence theorem gives $\varepsilon_1>0$ with $\liminf_{t\to\infty}\dist(x(t),\Om\cap\mathcal X_0)\ge\varepsilon_1$ for every $x_0\in X^{\circ}$. For $x=(z,y)\in\Om$ the point $(z,0)$ lies in $\Om\cap\mathcal X_0$, so $\dist(x,\Om\cap\mathcal X_0)\le\|y\|\le\one^{\!\top}y$, while $\dist(x,\Om\cap\mathcal X_0)\ge\dist(x,\mathcal X_0)=\|y\|\ge m^{-1/2}\one^{\!\top}y$. Hence \eqref{eq:persist} holds with $\varepsilon_0:=\varepsilon_1$. Ultimate boundedness is Theorem \ref{thm:region}, and the final statement is Theorem \ref{thm:EEexist}.
\end{proof}

\begin{corollary}[Extinction]\label{cor:extinct}
Under the hypotheses of Theorem \ref{thm:GASDFE} with $\Rz(\phi)<1$, the infection is eliminated at every exponential rate $\alpha<-s(F(\phi)-V(\phi))$, uniformly for initial data in $\Om$.
\end{corollary}

\begin{proof}
Immediate from \eqref{eq:expdecay}, whose constants depend only on $F(\phi)-V(\phi)$ and $\alpha$.
\end{proof}

Corollary \ref{cor:extinct} requires Condition \ref{cond:dom}. Without it, $\Rz(\phi)<1$ need not imply elimination, as Theorem \ref{thm:backward} shows by exhibiting endemic equilibria below threshold; the example of Section \ref{sec:numerics} has a locally stable one.

\section{Bifurcation analysis}\label{sec:bifurcation}

At $\Rz=1$ the disease-free equilibrium loses hyperbolicity and a branch of endemic equilibria emerges. Whether it is supercritical or subcritical is decided by the quadratic terms of the vector field. We carry out the center manifold reduction for the neural SEIR model \eqref{eq:seir} with $\beta_\phi\in C^2$ depending on all four compartments; the coefficients depend on the transition graph, so the computation is explicit, but the mechanism it exposes is not specific to \eqref{eq:seir}. A trained closure $\beta_\phi$ is a single function, not a family; to speak of a bifurcation we embed it in the one parameter family $\kappa\beta_\phi$, $\kappa>0$, which rescales the learned transmission rate uniformly and moves $\Rz$ through $1$ without retraining. Since \eqref{eq:R0seir} is linear in $\beta^\ast$, the family $\kappa\beta_\phi$ has reproduction number $\kappa\Rz(\phi)$, so the critical value is $\kappa^\ast=1/\Rz(\phi)$, written out in \eqref{eq:kstar}; the parameters $\Lambda,\mu,\sigma,\gamma,\delta$ and the weights $\phi$ are held fixed. Epidemiologically, $\kappa$ is a uniform proportional change in transmission, such as a reduction in contacts, and $\kappa=1$ is the trained model; $\kappa$ is a derived bifurcation parameter, not a quantity estimated from data. The statement below therefore concerns the family $\kappa\beta_\phi$, and the derivatives that enter it are those of the trained network at $E_0$. We use the reduction theorem of Castillo-Chavez and Song \cite{CastilloChavezSong2004}, a consequence of the center manifold theorem \cite{Carr1981,GuckenheimerHolmes1983}: for $\dot\chi=\mathcal H(\chi,\varphi)$ with $\mathcal H\in C^2$, $\mathcal H(0,\varphi)\equiv0$ and $A:=D_\chi\mathcal H(0,0)$ having a simple zero eigenvalue and all other eigenvalues with negative real part, let $\varpi$ and $\upsilon$ be right and left null vectors of $A$ with $\upsilon^{\!\top}\varpi=1$, where $\varpi_j>0$ in every coordinate $j$ in which the equilibrium of interest vanishes, while components in which it is positive may have either sign \cite[Theorem~4.1 and Remark~1]{CastilloChavezSong2004}, and put
\begin{equation}\label{eq:ab}
a:=\sum_{k,i,j}\upsilon_k\varpi_i\varpi_j\frac{\partial^2\mathcal H_k}{\partial\chi_i\partial\chi_j}(0,0),\qquad
b:=\sum_{k,i}\upsilon_k\varpi_i\frac{\partial^2\mathcal H_k}{\partial\chi_i\partial\varphi}(0,0);
\end{equation}
if $b>0$, the bifurcation at $\varphi=0$ is backward when $a>0$ and forward when $a<0$.

\begin{theorem}\label{thm:backward}
Consider \eqref{eq:seir} with $\beta_\phi\in C^2$ nonnegative and $\beta_\phi(E_0)>0$, replace $\beta_\phi$ by $\kappa\beta_\phi$ with $\kappa>0$, write $S^\ast=\Lambda/\mu$, $\beta^\ast:=\beta_\phi(E_0)$, and let $\beta_S,\beta_E,\beta_I,\beta_R$ be the partial derivatives of $\beta_\phi$ at $E_0$. Let
\begin{equation}\label{eq:kstar}
\kappa^{\ast}=\frac{\mu(\mu+\sigma)(\mu+\gamma+\delta)}{\sigma\beta^{\ast}\Lambda}
\end{equation}
be the value of $\kappa$ at which $\Rz=1$, and set $\varphi:=\kappa-\kappa^\ast$. With the null vectors in \eqref{eq:ab} normalized by $\varpi_3=\upsilon_3=1$, which differs from the normalization $\upsilon^{\!\top}\varpi=1$ by the common positive factor $\upsilon^{\!\top}\varpi$ in $a$ and $b$ and so leaves their signs unchanged, $b>0$ and
\begin{equation}\label{eq:acoef}
a=2\upsilon_2S^{\ast}\kappa^{\ast}\Bigl[\beta_I+\frac{\gamma}{\mu}\beta_R+\frac{S^{\ast}\beta^{\ast}\kappa^{\ast}}{\mu+\sigma}\beta_E-\frac{\kappa^{\ast}\beta^{\ast}}{\mu}\bigl(\beta^{\ast}+S^{\ast}\beta_S\bigr)\Bigr],\qquad \upsilon_2=\frac{\sigma}{\mu+\sigma}>0 .
\end{equation}
The bifurcation at $\Rz=1$ is therefore backward if
\begin{equation}\label{eq:backcrit}
\beta_I+\frac{\gamma}{\mu}\beta_R+\frac{S^{\ast}\beta^{\ast}\kappa^{\ast}}{\mu+\sigma}\beta_E\ >\ \frac{\kappa^{\ast}\beta^{\ast}}{\mu}\bigl(\beta^{\ast}+S^{\ast}\beta_S\bigr),
\end{equation}
and forward if the reverse strict inequality holds; in the degenerate case of equality $a=0$ and higher-order terms decide. If $\beta_\phi$ is constant then $a<0$ and the bifurcation is forward.
\end{theorem}

\begin{proof}
Put $\chi=(S-S^\ast,E,I,R)$, so the equilibrium is at the origin for every $\kappa$. The first derivatives of the incidence $\kappa\beta_\phi(x)SI$ are
\[
\partial_S=\kappa I(\beta_\phi+S\partial_S\beta_\phi),\quad \partial_E=\kappa SI\,\partial_E\beta_\phi,\quad \partial_I=\kappa S(\beta_\phi+I\partial_I\beta_\phi),\quad \partial_R=\kappa SI\,\partial_R\beta_\phi ;
\]
all except $\partial_I$ carry a factor $I$ and vanish at $E_0$, where $\partial_I=\kappa S^\ast\beta^\ast$. Hence at $\kappa=\kappa^\ast$
\[
J_0=\begin{pmatrix}
-\mu&0&-\kappa^\ast S^{\ast}\beta^{\ast}&0\\
0&-(\mu+\sigma)&\kappa^\ast S^{\ast}\beta^{\ast}&0\\
0&\sigma&-(\mu+\gamma+\delta)&0\\
0&0&\gamma&-\mu
\end{pmatrix},
\]
whence $\spec(J_0)=\{-\mu,-\mu\}\cup\spec(J^{(2)}_0)$ with $J^{(2)}_0=\begin{psmallmatrix}-(\mu+\sigma)&\kappa^\ast S^\ast\beta^\ast\\ \sigma&-(\mu+\gamma+\delta)\end{psmallmatrix}$. By \eqref{eq:kstar}, $\sigma\kappa^\ast S^\ast\beta^\ast=(\mu+\sigma)(\mu+\gamma+\delta)$, so $\det J_0^{(2)}=0$ and $\tr J^{(2)}_0<0$: the zero eigenvalue is simple and the remaining spectrum lies in the open left half-plane.

Solving $J_0\varpi=0$ with $\varpi_3=1$ gives $\varpi_1=-\kappa^\ast S^\ast\beta^\ast/\mu$, $\varpi_2=\kappa^\ast S^\ast\beta^\ast/(\mu+\sigma)$ and $\varpi_4=\gamma/\mu$, the third equation being satisfied identically by \eqref{eq:kstar}. The components $\varpi_2,\varpi_3,\varpi_4$ in the coordinates $E,I,R$, which vanish at $E_0$, are positive, and only $\varpi_1$, in the coordinate $S$ with $S^\ast>0$, is negative, as the reduction theorem permits. Solving $\upsilon^{\!\top}J_0=0$ gives $\upsilon_1=\upsilon_4=0$ from the first and fourth columns and $\upsilon_2=\sigma\upsilon_3/(\mu+\sigma)$ from the second; we normalize $\upsilon_3=1$, since rescaling to enforce $\upsilon^{\!\top}\varpi=1$ multiplies $a$ and $b$ by the same positive factor.

As $\upsilon_1=\upsilon_4=0$ and the third equation of \eqref{eq:seir} is linear, only the second contributes to $a$. With $\mathcal G(x):=\kappa^\ast\beta_\phi(x)SI$, differentiation at $I=0$ gives
\[
\frac{\partial^2\mathcal G}{\partial x_j\partial x_k}\Big|_{E_0}=\kappa^\ast\bigl(\delta_{jS}\delta_{kI}\beta^{\ast}+\delta_{kI}S^{\ast}\beta_{x_j}+\delta_{jI}\delta_{kS}\beta^{\ast}+\delta_{jI}S^{\ast}\beta_{x_k}\bigr),
\]
every second derivative of $\beta_\phi$ being multiplied by $I=0$. Contracting with $\varpi$ and using $\varpi_3=1$,
\[
\sum_{j,k}\varpi_j\varpi_k\frac{\partial^2\mathcal G}{\partial x_j\partial x_k}\Big|_{E_0}=2\kappa^\ast\Bigl(\varpi_1\beta^{\ast}+S^{\ast}\bigl(\beta_S\varpi_1+\beta_E\varpi_2+\beta_I+\beta_R\varpi_4\bigr)\Bigr),
\]
and multiplying by $\upsilon_2$ and substituting $\varpi$ gives \eqref{eq:acoef}. For $b$, the mixed derivatives $\partial^2\mathcal H_k/\partial\chi_i\partial\kappa$ are the $\chi_i$ derivatives of $\mp\beta_\phi SI$ in the $S$ and $E$ equations and vanish otherwise; by the computation above they vanish at $E_0$ unless $i$ is the $I$ coordinate, and since $\upsilon_1=0$ the only surviving mixed derivative weighted by a nonzero component of $\upsilon$ is $\partial^2(\kappa\beta_\phi SI)/\partial I\partial\kappa=S^\ast\beta^\ast$ in the second equation, so $b=\upsilon_2\varpi_3S^\ast\beta^\ast>0$. Criterion \eqref{eq:backcrit} is $a>0$, and for constant $\beta_\phi$ formula \eqref{eq:acoef} reduces to $a=-2\upsilon_2S^\ast\kappa^{\ast2}\beta^{\ast2}/\mu<0$.
\end{proof}

\begin{corollary}\label{cor:backprev}
If the learned transmission rate depends only on prevalence, $\beta_\phi=\beta_\phi(I)$, then \eqref{eq:seir} exhibits a backward bifurcation at $\Rz=1$ if
\begin{equation}\label{eq:backprev}
\frac{\beta'_\phi(0)}{\beta_\phi(0)}\ >\ \frac{(\mu+\sigma)(\mu+\gamma+\delta)}{\sigma\Lambda},
\end{equation}
where $\beta_\phi(0)$ and $\beta'_\phi(0)$ denote the value and the derivative of $I\mapsto\beta_\phi(I)$ at zero prevalence, that is, at $E_0$, and a forward bifurcation if the reverse strict inequality holds. In particular any nonincreasing dependence on prevalence yields a forward bifurcation.
\end{corollary}

\begin{proof}
Setting $\beta_S=\beta_E=\beta_R=0$ in \eqref{eq:backcrit} gives $\beta_I>\kappa^\ast\beta^{\ast2}/\mu$, and substituting \eqref{eq:kstar} with $\beta^\ast=\beta_\phi(0)$ and $S^\ast=\Lambda/\mu$ yields $\kappa^\ast\beta^{\ast2}/\mu=\beta_\phi(0)(\mu+\sigma)(\mu+\gamma+\delta)/(\sigma\Lambda)$.
\end{proof}

\begin{remark}\label{rem:consistency}
The same computation for the separable model \eqref{eq:sir} shows that its bifurcation at $\Rz=1$ is backward if $\beta'_\phi(0)/\beta_\phi(0)>(\mu+\gamma+\delta)/\Lambda$ and forward if the reverse strict inequality holds, the analogue of Corollary \ref{cor:backprev}. This is consistent with Theorem \ref{thm:GASEE}: in $\Ecal_{\mathrm{sep}}$ the ratio $h_\phi(I)/I=\beta_\phi(I)$ is nonincreasing, so $\beta'_\phi(0)\le0$, the bifurcation is forward, and the unique endemic equilibrium is globally stable. The sublinearity constraint therefore excludes bistability rather than merely simplifying the algebra. Conversely, a closure that has learned to amplify transmission at low prevalence, a plausible outcome when an accelerating early epidemic phase is attributed to transmission rather than to structured depletion of susceptibles, changes the bifurcation structure; Corollary \ref{cor:bifnonid} shows that this can happen without any change of the fitted trajectory, because the quantities entering \eqref{eq:backprev} are values and derivatives at $E_0$, which data supported away from $E_0$ do not constrain.
\end{remark}

\section{Sensitivity analysis}\label{sec:sensitivity}

\subsection{Normalized indices and neural elasticities}

Because $\Rz$ is available in closed form, its sensitivities are exact. Recall the normalized forward sensitivity index $\Upsilon_p:=(\partial\Rz/\partial p)(p/\Rz)=\partial\ln\Rz/\partial\ln p$ \cite{Chitnis2008}. For \eqref{eq:seir} the coefficient $\beta^\ast=\beta_\phi(\Lambda/\mu,0,0,0)$ depends on the demographic parameters through $S^\ast=\Lambda/\mu$, while the weights $\phi$, and hence the function $\beta_\phi$, are held fixed. We introduce the neural elasticity at the disease-free equilibrium,
\begin{equation}\label{eq:elasticity}
\varepsilon_S:=\frac{S^{\ast}\partial_S\beta_\phi(E_0)}{\beta_\phi(E_0)}=\frac{\partial\ln\beta_\phi}{\partial\ln S}\Big|_{E_0},
\end{equation}
assuming $\beta^\ast>0$. Taking logarithms of \eqref{eq:R0seir} gives
\[
\ln\Rz
=
\ln\sigma+\ln\beta^\ast+\ln\Lambda-\ln\mu
-\ln(\mu+\sigma)-\ln(\mu+\gamma+\delta).
\]
Since $S^\ast=\Lambda/\mu$ and $\beta^\ast=\beta_\phi(S^\ast,0,0,0)$, the chain rule and \eqref{eq:elasticity} give
\[
\frac{\partial\ln\beta^\ast}{\partial\ln\Lambda}
=\varepsilon_S,
\qquad
\frac{\partial\ln\beta^\ast}{\partial\ln\mu}
=-\varepsilon_S,
\]
because $\partial\ln S^\ast/\partial\ln\Lambda=1$ and $\partial\ln S^\ast/\partial\ln\mu=-1$; the parameters $\sigma,\gamma,\delta$ do not affect $S^\ast$ and enter only through the explicit factors in \eqref{eq:R0seir}. Therefore,
\begin{equation}\label{eq:indices}
\Upsilon_\Lambda=1+\varepsilon_S,\qquad \Upsilon_\mu=-(1+\varepsilon_S)-\frac{\mu}{\mu+\sigma}-\frac{\mu}{\mu+\gamma+\delta},\qquad \Upsilon_\sigma=\frac{\mu}{\mu+\sigma},
\end{equation}
\begin{equation}\label{eq:indices2}
\Upsilon_\gamma=-\frac{\gamma}{\mu+\gamma+\delta},\qquad \Upsilon_\delta=-\frac{\delta}{\mu+\gamma+\delta},\qquad \Upsilon_{\beta^{\ast}}=1 .
\end{equation}
Here $\Upsilon_{\beta^\ast}$ is the index with respect to a direct rescaling of the learned rate at $E_0$. The classical indices are recovered when $\varepsilon_S=0$. Each index is local: $\Upsilon_p$ is the approximate percentage change in $\Rz$ produced by a one percent change in $p$ at the nominal parameter values, the other parameters being fixed. The indices for $\sigma$, $\gamma$ and $\delta$ do not involve the closure, whereas those for $\Lambda$ and $\mu$ are shifted by $\varepsilon_S$. A ranking by $|\Upsilon_p|$ holds near the nominal values and need not persist over wide parameter ranges, which is the purpose of the global methods of Section \ref{sec:globalsens}. The neural elasticity therefore quantifies how state dependence learned by the closure modifies the classical parameter sensitivities. In particular, $\varepsilon_S>0$ amplifies the sensitivity of $\Rz$ to recruitment, whereas $\varepsilon_S<0$ attenuates it. If $\varepsilon_S<-1$, then $\Upsilon_\Lambda<0$, so increasing recruitment decreases $\Rz$, a behavior that cannot occur under constant mass-action transmission. Notably, the same threshold $\varepsilon_S=-1$ appears in the domination condition discussed after Condition \ref{cond:dom}; beyond this threshold the global stability guarantee of Theorem \ref{thm:GASDFE} may also fail.

\subsection{Sensitivity with respect to the weights}

\begin{proposition}\label{prop:neuralsens}
Under the hypotheses of Theorem \ref{thm:R0closed}, let $\phi\in\Phi\subset\R^d$, suppose that $\phi\mapsto\beta^\phi(x)$ and the progression rates are $C^1$ with derivatives continuous in $(x,\phi)$, and suppose that the uninfected part of the closure vanishes, so that $E_0$ does not depend on $\phi$ (Proposition \ref{prop:dfeindep}). First suppose that the learned closure affects transmission only, so that $V$ is independent of $\phi$. Then
\begin{equation}\label{eq:dR0dphi}
\frac{\partial\Rz}{\partial\phi_\ell} = z^{\ast}_{s_\ast} \Bigl( \frac{\partial\beta^{\phi}(E_0)}{\partial\phi_\ell} \Bigr)^{\!\top} V^{-1}w,
\end{equation}
which requires one gradient of the network output at $E_0$ and one linear solve.

If progression is also learned, then $V=V(\phi)$. Provided $\Rz(\phi)>0$,
\[
\frac{\partial\Rz}{\partial\phi_\ell}
=
\frac{\upsilon^{\!\top}
(\partial_{\phi_\ell}K)\omega}
{\upsilon^{\!\top}\omega},
\]
where
\[
\partial_{\phi_\ell}K
=
(\partial_{\phi_\ell}F)V^{-1}
-
FV^{-1}(\partial_{\phi_\ell}V)V^{-1},
\]
and $\omega=z^\ast_{s_\ast}w$ and $\upsilon=V^{-\!\top}\beta^\phi(E_0)$ are right and left eigenvectors of $K$ for the simple eigenvalue $\Rz$ (Remark \ref{rem:Vpert}); the quotient does not depend on their normalization.
\end{proposition}

\begin{proof}
The first formula follows from linearity of \eqref{eq:R0closed} in $\beta^\phi(E_0)$ for fixed $V$, the second from \eqref{eq:perron-deriv}, which applies because the nonzero eigenvalue of the rank-one matrix $K$ is simple (Remark \ref{rem:Vpert}), together with the chain rule and $\partial(V^{-1})=-V^{-1}(\partial V)V^{-1}$.
\end{proof}

If the uninfected part of the closure is also learned, then $E_0=E_0(\phi)$ and the derivatives above must be read as total derivatives, for instance $\frac{\dd}{\dd\phi_\ell}\beta^\phi(E_0(\phi))=\partial_{\phi_\ell}\beta^\phi(E_0)+D_x\beta^\phi(E_0)\,\partial_{\phi_\ell}E_0$, where $\partial_{\phi_\ell}z^\ast=-Dg_\phi(z^\ast)^{-1}\partial_{\phi_\ell}g_\phi(z^\ast)$ by differentiation of $g_\phi(z^\ast)=0$, which is legitimate because $Dg_\phi(z^\ast)$ is invertible by Assumption \ref{as:A5}; the factor $z^\ast_{s_\ast}$ in $F$ is then differentiated as well.

Formula \eqref{eq:dR0dphi} shows that the gradient of the threshold with respect to the weights is supported at the single point $E_0$. A loss depending only on the fitted trajectory therefore has, generically, no component along directions in weight space that alter $\Rz$ while preserving $\beta^\phi$ on the observed orbit; this is the differential counterpart of Theorem \ref{thm:nonid}.

\subsection{Global sensitivity}\label{sec:globalsens}

Latin hypercube sampling combined with partial rank correlation coefficients \cite{McKay1979,BlowerDowlatabadi1994,Marino2008} transfers to the present setting once the correct inputs are identified. By Theorem \ref{thm:R0closed}, $\Rz$ depends on the closure only through the finite dimensional vector $\beta^\ast:=\beta^\phi(E_0)$ and, when relevant, through the values of the progression closures at $E_0$. The design should therefore sample $\theta$ jointly with $\beta^\ast$ over ranges reflecting the posterior or bootstrap uncertainty of the trained model, rather than sampling the weights $\phi$, whose coordinates are neither interpretable nor independently constrained. Partial rank correlation measures monotone association, so it is appropriate only for outputs that are monotone in each input. Formula \eqref{eq:R0closed} is nondecreasing in each component of $\beta^\ast$, and for \eqref{eq:seir} with $\beta^\ast$ sampled as an independent input, \eqref{eq:indices}--\eqref{eq:indices2} with $\varepsilon_S=0$ show that $\Rz$ is also monotone in $\Lambda,\mu,\sigma,\gamma,\delta$. If instead $\beta^\ast$ is obtained by evaluating the trained network at $S^\ast=\Lambda/\mu$, the sign of $\Upsilon_\Lambda$ is that of $1+\varepsilon_S$, which may change over the sampled range, and in models with other structures $\Rz$ need not be monotone in every parameter. Monotonicity must then be checked, for instance with scatter plots, and variance-based measures used when it fails \cite{Marino2008}; the same applies to other outputs, such as peak prevalence. Three sources of uncertainty should be distinguished: uncertainty in the mechanistic parameters $\theta$; uncertainty in the network weights, which enters $\Rz$ only through $\beta^\ast$ and is propagated by sampling $\beta^\ast$ from its posterior or bootstrap distribution; and uncertainty in the extrapolation of the learned closure to states that the observed trajectory does not visit. The procedure addresses the first two. It does not address the third, which is the structural uncertainty analyzed next. Accordingly, the indices of this section describe the sensitivity of a given fitted model; they do not show that $\Rz$ or $\beta^\ast$ is determined by the data, which is the question of Section \ref{sec:identifiability}.

\section{Approximation and identifiability}\label{sec:identifiability}

Suppose the data are generated by \eqref{eq:ude} with an unknown closure $\beta^{\dagger}$ and let $\beta^\phi$ be a trained approximation.

\begin{theorem}[Error propagation]\label{thm:errorprop}
Let $x^{\dagger}$ and $x^{\phi}$ solve \eqref{eq:ude} on $[0,T]$ with values in $\Om$, with transmission closures $\beta^{\dagger}$ and $\beta^{\phi}$, all other terms of the vector field being identical, and with the same initial datum, and let $\sup_{x\in\Om}\|\beta^{\dagger}(x)-\beta^{\phi}(x)\|\le\varepsilon$. If $L_\Om$ is a Lipschitz constant on $\Om$ for the field with closure $\beta^\phi$, then
\begin{equation}\label{eq:errtraj}
\sup_{0\le t\le T}\|x^{\dagger}(t)-x^{\phi}(t)\|\ \le\ \sqrt2\Bigl(\frac{\Lambda}{\mu}\Bigr)^{2}\frac{\ee^{L_\Om T}-1}{L_\Om}\,\varepsilon,
\end{equation}
and, if the two models share the transfer matrix $V$, $|\Rz(\beta^{\dagger})-\Rz(\beta^{\phi})|\le\varepsilon(\Lambda/\mu)\one^{\!\top}V^{-1}w$. If the initial data differ, \eqref{eq:contdep} adds the term $\ee^{L_\Om T}\|x^\dagger(0)-x^\phi(0)\|$ to the right-hand side of \eqref{eq:errtraj}.
\end{theorem}

\begin{proof}
Write $F_\beta$ for the field \eqref{eq:ude} with transmission closure $\beta$. At any $x\in\Om$, the two models differ only through their transmission closures, so their vector fields differ by
\[
(\hat w-e_{s_\ast})x_{s_\ast}
\bigl(\beta^{\dagger}(x)-\beta^{\phi}(x)\bigr)^{\!\top}y.
\]
Since $\one^{\!\top}w=1$ and $w\ge0$, one has $\|\hat w-e_{s_\ast}\|\le\sqrt2$. Moreover, Theorem \ref{thm:region} gives $x_{s_\ast}\le\Lambda/\mu$ and $\|y\|\le\Lambda/\mu$. Therefore
\[
\|F_{\beta^\dagger}(x)-F_{\beta^\phi}(x)\|
\le
\sqrt2\,x_{s_\ast}
\|\beta^\dagger(x)-\beta^\phi(x)\|\,\|y\|
\le
\sqrt2\Bigl(\frac{\Lambda}{\mu}\Bigr)^2\varepsilon.
\]
Applying \eqref{eq:contdep} with equal initial data gives \eqref{eq:errtraj}. The second bound is \eqref{eq:R0lip} with $z^\ast_{s_\ast}\le\Lambda/\mu$.
\end{proof}

The bound on $\Rz$ is useful only if the closure is controlled at $E_0$. Trajectory data provide no such control. To make this precise, we fix the mechanistic model, $\theta$, $w$, and the progression and uninfected closures, so that $V$ is fixed and only the transmission closure varies. We write $\mathcal B$ for the resulting function class of all admissible transmission closures $\beta:\Rp^n\to[0,\bar G]^m$ satisfying the regularity requirements above. Thus $\mathcal B$ is a class of admissible functions, not a finite collection of fitted networks.

\begin{definition}[Identifiability from data on a set]\label{def:orbitid}
Let $\Gamma\subset\Om$ be compact. Two closures $\beta^{(1)},\beta^{(2)}\in\mathcal B$ are \emph{indistinguishable on $\Gamma$} if their corresponding vector fields agree at every state in $\Gamma$, that is,
\[
F_{\beta^{(1)}}(x)=F_{\beta^{(2)}}(x)
\qquad\text{for every }x\in\Gamma.
\]
A functional $Q:\mathcal B\to\R$ is \emph{identifiable from data on $\Gamma$} if
\[
Q(\beta^{(1)})=Q(\beta^{(2)})
\]
whenever $\beta^{(1)}$ and $\beta^{(2)}$ are indistinguishable on $\Gamma$.
\end{definition}

If $\Gamma$ is forward invariant, for instance the closure of a bounded orbit, then indistinguishable closures generate identical solutions from every initial condition in $\Gamma$, by Theorem \ref{thm:local}, and hence identical observations $Cx(t)$ for every observation operator $C$; the notion is therefore independent of which compartments are observed and of the sampling scheme. It describes what a single observed epidemic can determine. This differs from structural functional identifiability \cite{LomanBaker2025,LomanBrowningBaker2026}, which asks whether an unknown function is determined by the outputs generated from all admissible initial conditions: when trajectories from every initial condition in $\Om$ are available, and the full state is observed, the field $F_\beta$ is determined on $\Om$ from $\dot x$. This gives $\lambda=\beta^{\!\top}y$ wherever $x_{s_\ast}>0$, but not the individual coefficients $\beta_j$, since the representation \eqref{eq:foi} is not unique. The reproduction number is nevertheless identifiable in that setting, because \eqref{eq:foi} gives $\beta(z,0)=\nabla_y\lambda(z,0)$, so $\beta(E_0)$ is determined by $\lambda$ near $E_0$. The mechanistic parameters $\theta$ are held fixed here; their own structural identifiability, and the entanglement of $\beta^\phi(E_0)$ with $\theta$ visible in \eqref{eq:R0seir}, where only the product $\beta^\ast\Lambda/\mu$ enters, are questions of parametric identifiability that can be assessed with differential algebra tools \cite{Alsammani2026Tutorial}.

Three features of Definition \ref{def:orbitid} should be kept in mind. First, it compares vector fields, not coefficient functions: by \eqref{eq:foi}, $F_{\beta^{(1)}}=F_{\beta^{(2)}}$ on $\Gamma$ means $x_{s_\ast}(\beta^{(1)}-\beta^{(2)})^{\!\top}y=0$ there, that is, equality of the force of infection $\lambda=\beta^{\!\top}y$ wherever $x_{s_\ast}>0$. This does not imply $\beta^{(1)}=\beta^{(2)}$ on $\Gamma$ when two or more infected classes are infectious; it does when a single class $j$ carries the infection, as in \eqref{eq:seir}, and $x_{s_\ast}y_j>0$ on $\Gamma$. Second, it is a structural notion: the data are noise free and the vector field is known exactly on $\Gamma$, so it describes what remains undetermined even under ideal observation. Practical identifiability, which concerns estimation from finitely many noisy observations, is a different question and is not treated here (Section \ref{sec:discussion}). Third, since agreement of vector fields implies agreement of every observed output, nonidentifiability in the sense of Definition \ref{def:orbitid} holds a fortiori under partial observation. Conversely, the positive statements below, Theorem \ref{thm:quantid}, Corollary \ref{cor:identifiable} and Proposition \ref{prop:monoident}, assume that the data control $\beta$ itself on a set, which under partial observation requires a separate observability argument. The results of this section are therefore indistinguishability constructions and bounds, not a procedure for identifying $\beta$ from one trajectory.

\begin{theorem}[Nonidentifiability from trajectory data]\label{thm:nonid}
Let $\Gamma\subset\Om$ be compact with $E_0\notin\Gamma$ and $d:=\dist(E_0,\Gamma)>0$; for instance, $\Gamma$ may be the closure of a bounded orbit of \eqref{eq:ude} whose omega limit set does not contain $E_0$, or the range of a solution other than $E_0$ on a finite time interval. Let $\beta^{(1)}\in\mathcal B$ and let $\beta^{\circ}\in[0,\bar G]^m$ be any constant vector. Then there is $\beta^{(2)}\in\mathcal B$, of the same regularity class as $\beta^{(1)}$, such that
\begin{enumerate}[label=(\roman*),leftmargin=2.2em,itemsep=1pt]
\item $\beta^{(2)}=\beta^{(1)}$ outside the ball $B_{d/2}(E_0)$, in particular on $\Gamma$, so the two closures are indistinguishable on $\Gamma$; every solution of the first model whose range on $[0,T]$ lies in $\Gamma$ is also a solution of the second on $[0,T]$, and if $\Gamma$ is forward invariant, the two models generate identical solutions for every initial condition in $\Gamma$;
\item $\beta^{(2)}(E_0)=\beta^{\circ}$, so that $\Rz(\beta^{(2)})=z^\ast_{s_\ast}(\beta^{\circ})^{\!\top}V^{-1}w$ attains, as $\beta^{\circ}$ varies, every value in $[0,\,\bar G z^\ast_{s_\ast}\one^{\!\top}V^{-1}w]$.
\end{enumerate}
In particular, choosing $\beta^{\circ}$ with $z^\ast_{s_\ast}(\beta^{\circ})^{\!\top}V^{-1}w\neq\Rz(\beta^{(1)})$ shows that $\Rz$ is not identifiable from data on $\Gamma$. Both closures can be approximated to arbitrary accuracy $\varepsilon$ in $C^1(\Om)$ by admissible network closures $\varrho\circ\N$, with $\varrho$ the scaled logistic map of Section \ref{sec:model} applied componentwise and $\N$ a feedforward network with a smooth nonpolynomial activation. The approximants are admissible but no longer coincide exactly on $\Gamma$: for a solution of the first model with range in $\Gamma$ on $[0,T]$, the solutions of the two network models from the same initial datum differ by $O(\varepsilon)$ on $[0,T]$, while their reproduction numbers differ by at least $|\Rz(\beta^{(2)})-\Rz(\beta^{(1)})|-2\varepsilon z^\ast_{s_\ast}\one^{\!\top}V^{-1}w$. The phenomenon therefore persists within a neural class up to the approximation error, as approximate rather than exact indistinguishability.
\end{theorem}

\begin{proof}
Agreement on $\Gamma$ constrains $\beta$ only on $\Gamma$, and admissibility imposes only bounds and smoothness, so the values of $\beta$ on the ball $B_{d/2}(E_0)$, which is disjoint from $\Gamma$, are free. Concretely, let $\chi\in C^\infty(\R^n;[0,1])$ with $\chi(E_0)=1$ and $\chi=0$ outside $B_{d/2}(E_0)$, for instance $\chi(x)=\chi_0(2\|x-E_0\|/d)$ with a fixed $\chi_0\in C^\infty(\R;[0,1])$ such that $\chi_0=1$ near $0$ and $\chi_0=0$ on $[1,\infty)$, and set $\beta^{(2)}:=(1-\chi)\beta^{(1)}+\chi\,\beta^{\circ}$. Pointwise $\beta^{(2)}(x)$ is a convex combination of $\beta^{(1)}(x)\in[0,\bar G]^m$ and $\beta^{\circ}\in[0,\bar G]^m$, so $\beta^{(2)}\in\mathcal B$, and it has the regularity of $\beta^{(1)}$. Statement (i) is clear from the construction: a solution of the first model with range in $\Gamma$ on $[0,T]$ satisfies the second equation there, since the fields coincide on $\Gamma$, and Theorem \ref{thm:local} identifies it with the solution of the second model; if $\Gamma$ is forward invariant, every solution with initial condition in $\Gamma$ has this property. For (ii), $\chi(E_0)=1$ gives $\beta^{(2)}(E_0)=\beta^{\circ}$; Theorem \ref{thm:R0closed} gives the stated value of $\Rz(\beta^{(2)})$, which is linear in $\beta^{\circ}$ for the fixed $V$ and, by \eqref{eq:R0sandwich}, sweeps the stated interval as $\beta^{\circ}$ runs from $0$ to $\bar G\one$. For the final claim, the approximants must respect the constraint $0\le\beta\le\bar G$ on all of $\Rp^n$, which uniform approximation alone does not provide; we approximate through the output map. Fix $\eta\in(0,1)$ and put $\beta^{(k)}_\eta:=(1-\eta)\beta^{(k)}+\tfrac12\eta\bar G\one$, which takes values in $[\tfrac12\eta\bar G,(1-\tfrac12\eta)\bar G]^m$ and satisfies $\|\beta^{(k)}_\eta-\beta^{(k)}\|_{C^1(\Om)}\le c_k\eta$ with $c_k$ independent of $\eta$. The inverse $\varrho^{-1}(v)=\ln(v/(\bar G-v))$ is smooth on $(0,\bar G)$, so $g_k:=\varrho^{-1}\circ\beta^{(k)}_\eta$ is $C^2$ on a neighborhood of $\Om$, on which $\beta^{(k)}_\eta$ takes values in $(0,\bar G)^m$ by continuity, and universal approximation in $C^1$ on compact sets \cite{HornikStinchcombeWhite1990,Pinkus1999} gives networks $\N_k$ with $\|\N_k-g_k\|_{C^1(\Om)}\le\delta$. The closures $\tilde\beta^{(k)}:=\varrho\circ\N_k$ take values in $(0,\bar G)^m$ on all of $\R^n$ and are $C^2$, hence admissible, and since $\varrho$ and $\varrho'$ are Lipschitz, $\|\tilde\beta^{(k)}-\beta^{(k)}_\eta\|_{C^1(\Om)}\le C_k\delta$ with $C_k$ depending on $\bar G$ and $\|g_k\|_{C^1(\Om)}$. Choosing $\eta$ and then $\delta$ small gives $\|\tilde\beta^{(k)}-\beta^{(k)}\|_{C^1(\Om)}\le\varepsilon$.

Now let $x^{(1)}$ be a solution of the first exact model with range in $\Gamma$ on $[0,T]$, so that $\|x^{(1)}(t)-E_0\|\ge d$, and let $\tilde x^{(k)}$ be the solution of the $k$th network model with the same initial datum; it stays in $\Om$ by Theorem \ref{thm:region}. On the compact set $K:=\Om\setminus B_{d/2}(E_0)$ the exact closures coincide, so $\sup_K\|\tilde\beta^{(k)}-\beta^{(1)}\|\le\varepsilon$ for $k=1,2$, and as in the proof of Theorem \ref{thm:errorprop} the field of each network model differs from that of the first exact model by at most $\sqrt2(\Lambda/\mu)^2\varepsilon$ on $K$. Let $\tau_k$ be the supremum of the $t\le T$ with $\tilde x^{(k)}([0,t])\subset K$. On $[0,\tau_k]$, \eqref{eq:contdep} applied on $K$ gives $\|\tilde x^{(k)}-x^{(1)}\|\le C_T\varepsilon$ with $C_T:=\sqrt2(\Lambda/\mu)^2(\ee^{L_KT}-1)/L_K$, where $L_K$ is a Lipschitz constant on $K$ for the first exact model. If $C_T\varepsilon<d/2$, then $\tilde x^{(k)}(\tau_k)$ lies at distance greater than $d/2$ from $E_0$, and continuity forces $\tau_k=T$. Hence $\|\tilde x^{(1)}-\tilde x^{(2)}\|\le2C_T\varepsilon$ on $[0,T]$. The statement on the reproduction numbers is \eqref{eq:R0lip} applied twice.
\end{proof}

The defect is quantitative as well as qualitative: what controls the error in $\Rz$ is the regularity of the closure between the data and the disease-free state.

\begin{theorem}[Quantitative identifiability]\label{thm:quantid}
Let $\beta^{(1)},\beta^{(2)}\in\mathcal B$ be Lipschitz on $\Om$ with constant $L_\beta$, in the sense $\|\beta^{(k)}(x)-\beta^{(k)}(\tilde x)\|_\infty\le L_\beta\|x-\tilde x\|$, let $\Gamma\subset\Om$ be nonempty with $\sup_{x\in\Gamma}\|\beta^{(1)}(x)-\beta^{(2)}(x)\|_\infty\le\varepsilon$, and let $d:=\dist(E_0,\Gamma)$. Then
\begin{equation}\label{eq:quantid}
|\Rz(\beta^{(1)})-\Rz(\beta^{(2)})|\ \le\ z^\ast_{s_\ast}\one^{\!\top}V^{-1}w\,(\varepsilon+2L_\beta d).
\end{equation}
The term $L_\beta d$ cannot be removed: for the closures of Theorem \ref{thm:nonid} with $\beta^\circ\in\{0,\bar G\one\}$ one has $\varepsilon=0$, while one of the two choices gives $|\Rz(\beta^{(1)})-\Rz(\beta^{(2)})|\ge\tfrac12\bar Gz^\ast_{s_\ast}\one^{\!\top}V^{-1}w>0$. For these pairs \eqref{eq:quantid} therefore forces $L_\beta d\ge\bar G/4$, and the construction gives $L_\beta d\le\mathrm{Lip}(\beta^{(1)})\,d+c_0\bar G$ with a constant $c_0$ that depends only on the fixed cutoff profile. When $\mathrm{Lip}(\beta^{(1)})\,d\lesssim\bar G$, both sides of \eqref{eq:quantid} are thus of order $\bar Gz^\ast_{s_\ast}\one^{\!\top}V^{-1}w$, and the bound is sharp up to a constant factor.
\end{theorem}

\begin{proof}
For $\eta>0$ pick $x^\flat\in\Gamma$ with $\|x^\flat-E_0\|\le d+\eta$. Then
\begin{align*}
\|\beta^{(1)}(E_0)-\beta^{(2)}(E_0)\|_\infty&\le\|\beta^{(1)}(E_0)-\beta^{(1)}(x^\flat)\|_\infty+\|\beta^{(1)}(x^\flat)-\beta^{(2)}(x^\flat)\|_\infty\\
&\qquad+\|\beta^{(2)}(x^\flat)-\beta^{(2)}(E_0)\|_\infty\ \le\ 2L_\beta(d+\eta)+\varepsilon ,
\end{align*}
and \eqref{eq:R0lip} gives \eqref{eq:quantid} after letting $\eta\downarrow0$. For the last statement, the bump in the proof of Theorem \ref{thm:nonid} is $\chi(x)=\chi_0(2\|x-E_0\|/d)$, whose Lipschitz constant is $2\|\chi_0'\|_\infty/d$. Writing $\beta^{(2)}(x)-\beta^{(2)}(\tilde x)=(1-\chi(x))(\beta^{(1)}(x)-\beta^{(1)}(\tilde x))+(\chi(x)-\chi(\tilde x))(\beta^\circ-\beta^{(1)}(\tilde x))$ and using $\|\beta^\circ-\beta^{(1)}\|_\infty\le\bar G$ gives $\mathrm{Lip}(\beta^{(2)})\le\mathrm{Lip}(\beta^{(1)})+2\|\chi_0'\|_\infty\bar G/d$, so $c_0=2\|\chi_0'\|_\infty$. Since $V^{-1}w\ge0$, the reproduction numbers for $\beta^\circ=0$ and $\beta^\circ=\bar G\one$ differ from $\Rz(\beta^{(1)})$ by $z^\ast_{s_\ast}\beta^{(1)}(E_0)^{\!\top}V^{-1}w$ and $z^\ast_{s_\ast}(\bar G\one-\beta^{(1)}(E_0))^{\!\top}V^{-1}w$, whose sum is $\bar Gz^\ast_{s_\ast}\one^{\!\top}V^{-1}w$; one of them is at least half of it, and this quantity is positive because $V^{-1}$ has positive diagonal (Proposition \ref{prop:ngm}) and $w\gneq0$.
\end{proof}

If $\beta^{(1)}$ and $\beta^{(2)}$ depend on $x$ only through a set of coordinates $x_J=(x_j)_{j\in J}$ and are Lipschitz in $x_J$ with constant $L_\beta$, the same proof gives \eqref{eq:quantid} with $d$ replaced by $d_J:=\dist((E_0)_J,\{x_J:x\in\Gamma\})\le d$. For prevalence-dependent closures of \eqref{eq:seir} this is $d_I=\inf_{x\in\Gamma}I$, so observations at low prevalence suffice to make it small; in general, $d_J$ is small only if the data come close to $E_0$ in every coordinate on which the closure depends, and $d$ only if they come close in all coordinates.

The hypothesis of Theorem \ref{thm:quantid} bounds the difference of the coefficient vectors on $\Gamma$, not merely the difference of the fluxes. By the discussion after Definition \ref{def:orbitid}, the two coincide when a single infected class carries the infection and $x_{s_\ast}y_j>0$ on $\Gamma$, but in general agreement of trajectories does not supply the componentwise bound, which is then an additional assumption on the hypothesis class. The theorem also presupposes, as everywhere in this section, that the two models share $z^\ast$ and $V$.

\begin{remark}\label{rem:lipcert}
For closures $\beta^\phi=\varrho\circ\N_\phi$, Lemma \ref{lem:nnlip} bounds $L_\beta$ by $\mathrm{Lip}(\varrho)L_\varsigma^{L-1}\prod_\ell\|W_\ell\|$, and sharper estimates of the Lipschitz constant of a trained network are available \cite{Virmaux2018}. A spectral norm penalty or a Lipschitz constraint during training therefore yields, through \eqref{eq:quantid}, a computable certificate on the reproduction number: two fits that agree to within $\varepsilon$ on the data cannot disagree on $\Rz$ by more than $z^\ast_{s_\ast}\one^{\!\top}V^{-1}w(\varepsilon+2L_\beta d)$. Regularization, priors and cross-validation do not create information about the vector field in a region the trajectory never enters, but they restrict the hypothesis class, and \eqref{eq:quantid} measures exactly what the restriction buys: the certificate is informative only when $L_\beta d$ is small compared with $\beta^\phi(E_0)$, that is when the closure is smooth on the scale of the distance between the data and the disease-free state.
\end{remark}

\begin{corollary}\label{cor:identifiable}
Suppose every closure in $\mathcal B$ consistent with the data satisfies $\sup_{D}\|\beta-\beta^\phi\|_\infty\le\varepsilon$ for some set $D\subset\Om$ and some reference closure $\beta^\phi$. If $E_0\in D$, all such closures give reproduction numbers within $2\varepsilon z^\ast_{s_\ast}\one^{\!\top}V^{-1}w$ of one another. If instead $\operatorname{dist}(E_0,D)>0$ and consistency with the data constrains the closure only on $D$, the data are compatible with every value of $\Rz$ in $[0,\bar Gz^\ast_{s_\ast}\one^{\!\top}V^{-1}w]$.
\end{corollary}

\begin{proof}
The first claim is \eqref{eq:R0lip} applied to two such closures, which differ by at most $2\varepsilon$ at $E_0$. For the second, apply the construction of Theorem \ref{thm:nonid} with a bump supported in a ball around $E_0$ disjoint from $\overline D$ and with $\beta^{\circ}=c\bar G\one$, $c\in[0,1]$; by Theorem \ref{thm:R0closed} the resulting reproduction number is $c\,\bar Gz^\ast_{s_\ast}\one^{\!\top}V^{-1}w$, which sweeps the stated interval.
\end{proof}

The same construction shows that the type of the bifurcation at threshold is not identifiable either, because the criterion \eqref{eq:backprev} involves the value and the first derivative of the learned rate at $E_0$.

\begin{corollary}[Nonidentifiability of the bifurcation type]\label{cor:bifnonid}
Consider \eqref{eq:seir} with a prevalence-dependent closure $\beta_\phi=\beta_\phi(I)$, let $\Gamma\subset\Om$ be compact with $d_I:=\inf_{x\in\Gamma}I>0$, let $\beta^{(1)}=\beta^{(1)}(I)$ be admissible, of class $C^2$, with bound $\bar G$, and let $b_0\in(0,\bar G)$ and $b_1\in\R$. Then there is an admissible $\beta^{(2)}=\beta^{(2)}(I)$ of class $C^2$ with bound $\bar G$ which coincides with $\beta^{(1)}$ for $I\ge d_I/2$, hence is indistinguishable from $\beta^{(1)}$ on $\Gamma$, and satisfies $\beta^{(2)}(0)=b_0$ and $(\beta^{(2)})'(0)=b_1$. Consequently, from data on $\Gamma$, the reproduction number \eqref{eq:R0seir} and the sign of the bifurcation coefficient $a$ in \eqref{eq:acoef} can be prescribed independently of one another: two closures indistinguishable on $\Gamma$ may share the same $\Rz$ and yet undergo a forward and a backward bifurcation, respectively.
\end{corollary}

\begin{proof}
Choose $r\in(0,d_I/2)$ with $r|b_1|\le\tfrac12\min\{b_0,\bar G-b_0\}$, put $p(I):=b_0+b_1I$, so that $p(I)\in[\tfrac12b_0,\tfrac12(\bar G+b_0)]\subset(0,\bar G)$ for $|I|\le r$, and let $\chi\in C^\infty(\R;[0,1])$ be even with $\chi=1$ on $[-r/2,r/2]$ and $\chi=0$ outside $(-r,r)$. Define $\beta^{(2)}(I):=(1-\chi(I))\beta^{(1)}(I)+\chi(I)p(I)$. For $I\ge0$ this is a convex combination of values in $[0,\bar G]$, so $\beta^{(2)}$ is admissible with bound $\bar G$; it is $C^2$, equals $\beta^{(1)}$ for $I\ge r$, and equals $p$ on $[0,r/2]$, whence $\beta^{(2)}(0)=b_0$ and $(\beta^{(2)})'(0)=b_1$. By \eqref{eq:R0seir}, $\Rz$ depends on the closure only through $b_0$, and by Corollary \ref{cor:backprev} the bifurcation is backward if $b_1/b_0>(\mu+\sigma)(\mu+\gamma+\delta)/(\sigma\Lambda)$ and forward if the reverse strict inequality holds; both can be set at will.
\end{proof}

Corollary \ref{cor:bifnonid} is a statement about the admissible function class. Within a neural class, the $C^1$ approximation in Theorem \ref{thm:nonid} controls the value and the first derivatives of the closure at $E_0$, on which the coefficient $a$ of Theorem \ref{thm:backward} depends continuously when $\beta^\ast>0$; a strict sign of $a$, and with it the bifurcation type, is therefore preserved by sufficiently accurate network approximants.

Monotone architectures restore identifiable information about the threshold, but only in one direction, and the direction is dictated by the monotonicity.

\begin{proposition}[Monotone bounds]\label{prop:monoident}
Let $\Ucal_+:=\{j\in\Ucal:z^\ast_j>0\}$ and $\Ucal_0:=\Ucal\setminus\Ucal_+$ be the uninfected classes that are occupied, respectively empty, at $E_0$; for \eqref{eq:seir}, $\Ucal_+=\{S\}$ and $\Ucal_0=\{R\}$. Let $x^\flat\in\Om$ with $x^\flat_j\le z^\ast_j$ for $j\in\Ucal_+$, and suppose the data determine the coefficient vector at $x^\flat$ to within $\varepsilon$, in the following sense: there is a vector $\hat\beta\in\R^m$, computed from the data alone, such that $\|\beta^\phi(x^\flat)-\hat\beta\|_\infty\le\varepsilon$ for every closure $\beta^\phi$ in the hypothesis class that is consistent with the data.
\begin{enumerate}[label=(\roman*),leftmargin=2.2em,itemsep=1pt]
\item If $\beta^\phi$ is componentwise nonincreasing in each coordinate of $\Ical\cup\Ucal_0$ and nondecreasing in each coordinate of $\Ucal_+$ on $\Om$, then
\begin{equation}\label{eq:monolower}
\Rz(\phi)\ \ge\ z^\ast_{s_\ast}\bigl(\hat\beta-\varepsilon\one\bigr)_+^{\!\top}V^{-1}w
\end{equation}
is an identifiable lower bound, where $(\cdot)_+$ denotes the componentwise positive part.
\item If $\beta^\phi$ is componentwise nondecreasing in each coordinate of $\Ical\cup\Ucal_0$ and nonincreasing in each coordinate of $\Ucal_+$ on $\Om$, then $\Rz(\phi)\le z^\ast_{s_\ast}(\hat\beta+\varepsilon\one)^{\!\top}V^{-1}w$ is an identifiable upper bound.
\end{enumerate}
\end{proposition}

The bounds are identifiable in the sense that their values depend only on $\hat\beta$, $\varepsilon$, $z^\ast$ and $V$, and so are the same for every closure in the class that is consistent with the data. The hypothesis on $\hat\beta$ is an observational one. For \eqref{eq:seir}, where a single class carries the infection, $\beta_\phi(x^\flat)$ equals the incidence $\beta_\phi SI$ divided by $S^\flat I^\flat$, so it is determined at any state with $S^\flat I^\flat>0$ at which the incidence and the state are observed, with $\varepsilon$ accounting for observation error. When several infected classes are infectious, the data determine only $(\beta^\phi)^{\!\top}y$ at $x^\flat$, and a componentwise $\hat\beta$ requires additional information (see the discussion after Definition \ref{def:orbitid}).

\begin{proof}
Let $\bar x$ be obtained from $x^\flat$ by setting the coordinates in $\Ical\cup\Ucal_0$ to zero; then $\bar x\in\Om$ and $\bar x\le x^\flat$ in these coordinates, while $\bar x$ and $x^\flat$ agree on $\Ucal_+$. In case (i), monotonicity in the coordinates of $\Ical\cup\Ucal_0$ gives $\beta^\phi(x^\flat)\le\beta^\phi(\bar x)$, and monotonicity in the coordinates of $\Ucal_+$, along which $\bar x\le E_0$ with equality elsewhere, gives $\beta^\phi(\bar x)\le\beta^\phi(E_0)$; hence $\beta^\phi(E_0)\ge\beta^\phi(x^\flat)\ge\hat\beta-\varepsilon\one$ by the hypothesis on $\hat\beta$, and also $\beta^\phi(E_0)\ge0$. The paths used stay in $\Om$, because decreasing coordinates preserves membership in $\Om$ and every point between $\bar x$ and $E_0$ is componentwise at most $E_0\in\Om$. The lower bound in \eqref{eq:R0sandwich} with $\beta^-=(\hat\beta-\varepsilon\one)_+$ gives \eqref{eq:monolower}. Case (ii) is symmetric, using the upper bound in \eqref{eq:R0sandwich}.
\end{proof}

Case (i) is the monotonicity that appears in the stability theory: nonincreasing dependence on the infected coordinates is a hypothesis of Theorem \ref{thm:GASEE}, gives the forward bifurcation in the last statement of Corollary \ref{cor:backprev}, and is part of the sufficient condition for the first inequality of \eqref{eq:dom} given after Condition \ref{cond:dom}. For \eqref{eq:seir}, an architecture of type (i), nonincreasing in $E$, $I$, $R$ and nondecreasing in $S$, bifurcates forward by \eqref{eq:backcrit}, satisfies the first inequality of \eqref{eq:dom}, and certifies the threshold from below, which is the conservative direction for control: the data prove that $\Rz$ is at least \eqref{eq:monolower}. An upper bound, case (ii), requires the opposite monotonicity, transmission amplified at higher prevalence and damped at higher susceptible density, and is not available from the architectures that guarantee the classical qualitative behavior. Data at the disease-free equilibrium itself give both bounds at once, by \eqref{eq:R0lip}.

Theorems \ref{thm:nonid} and \ref{thm:quantid} concern the geometry of the inverse problem, not a particular estimation procedure. A threshold reported from a fitted UDE should therefore be accompanied by a statement of what constrains the closure near the disease-free state: observations close to the disease-free state, which make $d$ small in \eqref{eq:quantid} (for prevalence-dependent closures, observations at low prevalence suffice, by the coordinate version of \eqref{eq:quantid} given after Theorem \ref{thm:quantid}); a Lipschitz bound on the closure, which makes $L_\beta d$ small; a mechanistic prior on transmission at low incidence; or a monotone architecture, which by Proposition \ref{prop:monoident} converts observations at higher prevalence into a one-sided bound at $E_0$.

\section{A numerical illustration}\label{sec:numerics}

The constructions of Theorem \ref{thm:nonid} and Corollary \ref{cor:bifnonid} are explicit, so their consequences can be computed exactly rather than inferred from a training run. We take the neural SEIR model \eqref{eq:seir} with a prevalence-dependent closure $\beta_\phi=\beta_\phi(I)$ and parameters $\Lambda=1$, $\mu=0.01$, $\sigma=0.2$, $\gamma=0.1$, and $\delta=0$ (time in days), so $S^\ast=100$, with bound $\bar G=0.01$.

Because the closure depends only on prevalence $I$, evaluation at the disease-free equilibrium $E_0=(S^\ast,0,0,0)$ reduces to evaluation at $I=0$. Thus
\[
\beta_\phi(E_0)=\beta_\phi(0),
\qquad
\partial_I\beta_\phi(E_0)=\beta_\phi'(0).
\]
Consequently,
\[
\Rz=K_0\,\beta_\phi(0),
\qquad
K_0=
\frac{\sigma S^\ast}
{(\mu+\sigma)(\mu+\gamma+\delta)}
=865.8,
\]
and, by Corollary \ref{cor:backprev}, the bifurcation is backward if
\[
\frac{\beta_\phi'(0)}{\beta_\phi(0)}
>
\frac{(\mu+\sigma)(\mu+\gamma+\delta)}
{\sigma\Lambda}
=0.1155,
\]
and forward if the reverse strict inequality holds. Three prevalence-dependent closures are compared. They are constructed, not fitted: closure 1 is a standard saturating incidence, and closures 2 and 3 are obtained from it by the cutoff construction of Corollary \ref{cor:bifnonid}, so that all three coincide exactly at high prevalence, while closure 2 changes the value of the rate at $E_0$, and with it the threshold (Theorem \ref{thm:nonid}; its slope at $E_0$ changes as well), and closure 3 changes only the slope at $E_0$, isolating the bifurcation type (Corollary \ref{cor:bifnonid}).

\begin{enumerate}[label=(\alph*),leftmargin=2.2em,itemsep=1pt]
\item Closure 1, the reference: saturating incidence $\beta^{(1)}(I)=b_0/(1+\alpha I)$ with $b_0=2.3\times10^{-3}$, $\alpha=0.05$. 
It is positive and nonincreasing, so Condition~\ref{cond:sublin}
holds (see Remark~\ref{rem:A6examples}), and the bifurcation
is forward by Corollary~\ref{cor:backprev}.
\item Closure 2: $\beta^{(2)}=(1-\chi)\beta^{(1)}+\chi\,p$ as in the proof of Corollary \ref{cor:bifnonid}, with $p(I)=b_0'+b_1'I$, $b_0'=9.24\times10^{-4}$ and $b_1'=2.5\times10^{-4}$.
\item Closure 3: the same construction with $p(I)=b_0+b_1I$, $b_1=5\times10^{-4}$, so that $\beta^{(3)}(0)=\beta^{(1)}(0)$.
\end{enumerate}
The reference trajectory starts at $x_0=(S,E,I,R)=(90,2,3,5)$ and is followed for $T=400$ days; it converges to the endemic equilibrium of closure 1, with $I^\ast=3.541$. Along this trajectory, the prevalence stays in $[3.000,\,6.731]$ and the Euclidean distance to $E_0$ is at least $11.7$. The bump $\chi$ is a smooth function of $I$ alone, equal to $1$ on $[0,r/2]$ and to $0$ for $I\ge r$, namely $\chi(I)=1-g(2|I|/r-1)$ with $g(t)=\psi(t)/(\psi(t)+\psi(1-t))$, $\psi(t)=\ee^{-1/t}$ for $t>0$ and $\psi(t)=0$ for $t\le0$, and $r=1.35<\tfrac12\min_\Gamma I$, so both modified closures are admissible with bound $\bar G$ and coincide exactly with closure 1 for $I\ge r=1.35$, a range that contains the visited prevalence range $[3.000,\,6.731]$. All three models were integrated from $x_0$ with an explicit Runge--Kutta method of order 8 (DOP853, relative tolerance $10^{-10}$, absolute tolerance $10^{-12}$; SciPy 1.18.1, NumPy 2.5.3, Python 3.14); the reported values are unchanged at relative tolerances $10^{-8}$ and $10^{-12}$ and with the implicit Radau and LSODA integrators. Because the three closures coincide on the visited prevalence range, uniqueness implies that the three models generate the same trajectory from $x_0$; the three numerical solutions are in fact identical, since the closures are evaluated only where they coincide. Endemic equilibria were computed from the scalar equation $\Lambda-\mu S(I)-cI=0$, with $c=(\mu+\sigma)(\mu+\gamma+\delta)/\sigma$ and $S(I)=c/(\kappa\beta(I))$, obtained by eliminating $S$, $E$ and $R$ from the equilibrium equations, and their stability from the eigenvalues of the Jacobian of \eqref{eq:seir}; the coefficients $a$ were also confirmed by a direct numerical evaluation of \eqref{eq:ab}. The quantities that the theory declares unconstrained by the trajectory differ as follows. The reason is that the trajectory samples the closures only for $I\in[3.000,\,6.731]$, whereas $\Rz$ depends on $\beta(0)$ alone (Theorem \ref{thm:R0closed}) and the sign of $a$ on $\beta(0)$ and $\beta'(0)$ (Corollary \ref{cor:backprev}), which are fixed by the closure on $[0,r/2]$, a range the trajectory never enters.

\begin{table}[t]
\centering
\caption{Three closures for the neural SEIR model \eqref{eq:seir} with $\Lambda=1$, $\mu=0.01$, $\sigma=0.2$, $\gamma=0.1$, $\delta=0$, which generate the same trajectory from $x_0=(90,2,3,5)$. Here $\Rz=K_0\beta(0)$ with $K_0=865.8$ (Theorem \ref{thm:R0closed}); the bifurcation of the family $\kappa\beta$ at $\kappa\Rz=1$ is backward if $\beta'(0)/\beta(0)>0.1155$ and forward if $\beta'(0)/\beta(0)<0.1155$ (Corollary \ref{cor:backprev}), equivalently $a>0$ or $a<0$, where $a$ is \eqref{eq:acoef} with $\beta_S=\beta_E=\beta_R=0$ and the normalization $\varpi_3=\upsilon_3=1$. The last column lists the equilibria of the closures as given ($\kappa=1$): the disease-free equilibrium $E_0$ and the endemic equilibria, identified by their prevalence $I^\ast$, with local stability determined from the eigenvalues of the Jacobian.}
\label{tab:numerics}
\footnotesize
\setlength{\tabcolsep}{3.5pt}
\begin{tabular}{@{}cccccccl@{}}
\toprule
Closure & $\beta(0)$ & $\beta'(0)$ & $\beta'(0)/\beta(0)$ & $\Rz$ & $a$ & Bifurcation & Equilibria at $\kappa=1$ \\
 & $(\times10^{-3})$ & $(\times10^{-4})$ & & & $(\times10^{-2})$ & & \\
\midrule
1 & 2.300 & $-$1.15 & $-$0.050 & 1.991 & $-$3.64 & forward & $E_0$ unstable; $3.541$ stable \\
2 & 0.924 & \phantom{$-$}2.50 & \phantom{$-$}0.271 & 0.800 & \phantom{$-$}3.41 & backward & $E_0$ stable; $0.879$ unstable; $3.541$ stable \\
3 & 2.300 & \phantom{$-$}5.00 & \phantom{$-$}0.217 & 1.991 & \phantom{$-$}2.24 & backward & $E_0$ unstable; $3.541$ stable \\
\bottomrule
\end{tabular}
\end{table}

\begin{figure}[t]
\centering
\includegraphics[width=\textwidth]{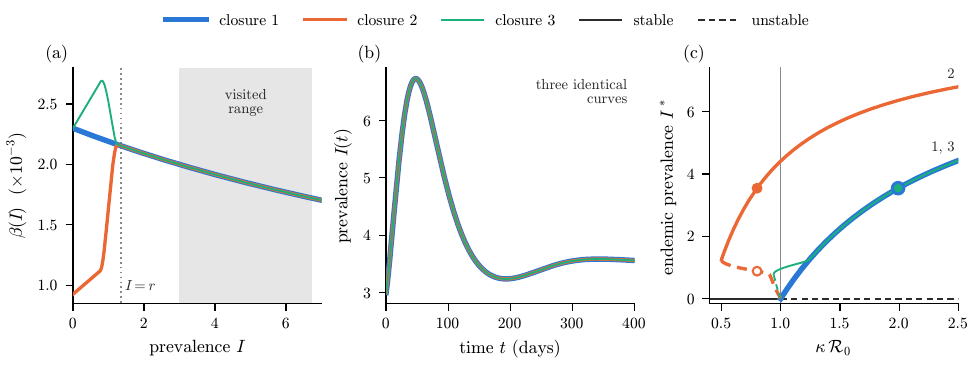}
\caption{(a) The three transmission closures as functions of prevalence; the shaded band is the prevalence range $[3.000,\,6.731]$ visited by the trajectory, and the dotted line marks $I=r=1.35$, beyond which the closures coincide exactly. (b) Prevalence $I(t)$ of the three models from the common initial datum $x_0=(90,2,3,5)$; the three solutions are identical and are drawn with decreasing line widths so that all remain visible. (c) Endemic equilibria $I^\ast$ of the rescaled families $\kappa\beta^{(k)}$ against their reproduction number $\kappa\Rz(\beta^{(k)})$, together with the disease-free branch $I^\ast=0$; solid branches are locally stable and dashed branches unstable, as determined by the eigenvalues of the Jacobian of \eqref{eq:seir}. Markers show the equilibria at $\kappa=1$ listed in Table \ref{tab:numerics} (filled: stable; open: unstable). Closure 1 bifurcates forward; closures 2 and 3 bifurcate backward, with turning points at $\kappa\Rz=0.50$ and $0.94$, and possess sub-threshold endemic equilibria. Each closure has the same color in all panels.}
\label{fig:nonid}
\end{figure}

Closures 1 and 2 illustrate Theorem \ref{thm:nonid}: the same trajectory is generated by a model with $\Rz=1.99$ and by a model with $\Rz=0.80$. The second model has a locally stable disease-free equilibrium by Theorem \ref{thm:LASDFE} and nevertheless converges, from $x_0$, to the endemic state $I^\ast=3.541$, which it shares with closure 1 because the two fields coincide there; an unstable endemic equilibrium at $I^\ast=0.879$, a saddle with a single unstable eigenvalue, separates the two attractors. Its one-dimensional unstable manifold connects it to $E_0$ in one direction and to the endemic state in the other, so it lies on the boundary of both basins of attraction. Theorem \ref{thm:backward} gives, for $a>0$, only the unstable endemic equilibrium just below threshold; the stable upper branch, and with it the bistability, is established here numerically: Figure \ref{fig:nonid}(c) shows that the sub-threshold branch of closure 2 persists down to $\kappa\Rz=0.50$, where it turns back. Closures 1 and 3 illustrate Corollary \ref{cor:bifnonid}: they share the trajectory and the reproduction number, and at $E_0$ they differ only in the derivative of the rate, which reverses the sign of $a$ and creates a sub-threshold endemic branch for $0.94<\kappa\Rz<1$ (Figure \ref{fig:nonid}(c)). None of this involves optimization: the three closures are equally consistent with the data by construction, and no loss function evaluated on the trajectory can distinguish them. The argument behind \eqref{eq:quantid} explains why the modifications must be steep. Since the closures depend only on prevalence, its coordinate version applies with $d_I=3.000$: closure 2 differs from closure 1 by $1.38\times10^{-3}$ at $I=0$ and coincides with it on $\Gamma$, so the two cannot share a Lipschitz constant below $1.376\times10^{-3}/(2\cdot3.000)=2.29\times10^{-4}$, about twice the maximal slope $1.15\times10^{-4}$ of closure 1 on $[0,1.35]$. Closure 3 has the same value at $E_0$ as closure 1, so \eqref{eq:quantid} imposes nothing on it; its steepness is forced by the prescribed slope at $E_0$. The maximal slopes $|\beta'|$ of closures 2 and 3 on $[0,1.35]$ are $3.1\times10^{-3}$ and $1.8\times10^{-3}$, against $8.7\times10^{-5}$ for closure 1 on the visited range. A Lipschitz constraint at the level $1.15\times10^{-4}$ of the reference closure's own regularity would exclude closures 2 and 3 while leaving closure 1 unaffected, which is the mechanism behind Remark \ref{rem:lipcert}. Closures 2 and 3 are mathematical possibilities, not estimated mechanisms: no biological interpretation is claimed for their behavior at low prevalence, and the example contains no epidemic data. Whether trained networks in practice land on closures of the second or third kind is an empirical question that the present construction does not answer; it shows only that nothing in the data forbids them.

\section{Discussion}\label{sec:discussion}

\paragraph{Architecture.} The admissible class turns two qualitative requirements into architectural facts. Flux parameterization with a factor $x_{s_k}$ on every outflow makes positivity and the population balance inequality \eqref{eq:Nbound} hold for every weight vector, so a model in the class is an epidemic model at initialization and stays one during training; the factorized force of infection \eqref{eq:foi} makes the disease-free set invariant and gives the model a next-generation matrix and a reproduction number for every weight vector. Neither property is obtained by penalties, and both survive retraining. The remaining classical hypotheses, domination \eqref{eq:dom}, sublinearity, and monotonicity, are conditions on the compact set $\Om$; monotone and sublinear architectures enforce them by construction, and for other architectures they are properties of the trained network that can in principle be certified, although certification of monotonicity or Lipschitz constants for deep networks is computationally demanding \cite{Virmaux2018}.

\paragraph{Interpretation of fitted closures.} The theory separates what a fitted closure determines on the region visited by the data from what it asserts elsewhere. Post-invasion dynamics, the endemic level, the shape of the epidemic curve, are governed by the closure on the visited region and are constrained by the data. The invasion threshold and the direction of the bifurcation at threshold are governed by the value and, as shown for \eqref{eq:seir}, the first derivatives of the closure at $E_0$ (Theorems \ref{thm:R0closed} and \ref{thm:backward}), and Theorem \ref{thm:nonid}, Corollary \ref{cor:bifnonid} and Section \ref{sec:numerics} show that trajectory data supported away from $E_0$ leave them free. A reproduction number or a bifurcation type read off a trained UDE is therefore a statement about the architecture and the regularization as much as about the data, and it should be reported together with what constrains the closure near the disease-free state.

\paragraph{Inference of $\Rz$ and data design.} Theorem \ref{thm:quantid} is the practical form of this statement. Two fits that agree to within $\varepsilon$ on the data cannot disagree on $\Rz$ by more than $z^\ast_{s_\ast}\one^{\!\top}V^{-1}w(\varepsilon+2L_\beta d)$. The three quantities are under the modeler's control to different degrees: $\varepsilon$ through the fit, $L_\beta$ through a Lipschitz or spectral norm constraint on the network, and $d$, or its coordinate version $d_J$, through the data, which is small only if the observations come close to the disease-free state in the coordinates on which the closure depends: for a prevalence-dependent closure a low-prevalence phase suffices, whereas a closure that also depends on the susceptible density requires data with $S$ near $S^\ast$, as at the beginning of an outbreak. Mechanistic priors on transmission at low incidence act on the same term. Monotone architectures of type (i) in Proposition \ref{prop:monoident}, which for \eqref{eq:seir} also guarantee a forward bifurcation, certify $\Rz$ from below; the upper bound requires the opposite monotonicity and is not available from those architectures. These are statements about what the data and the hypothesis class determine, not about a particular estimator, and they apply equally to point estimates, Bayesian posteriors, and cross-validated fits, none of which can identify the vector field outside the data-supported region without additional structural assumptions, regularization, or prior information.

\paragraph{Interventions.} The direction of the bifurcation decides whether reducing $\Rz$ below one suffices for elimination or whether a sub-threshold endemic state persists. Two learned models that fit the same data and report the same $\Rz$ may differ in this respect (Corollary \ref{cor:bifnonid}, Table \ref{tab:numerics}). The mathematical results do not say which model is right; they say that the data do not decide, and they identify the quantity, the derivative of the learned rate at $E_0$, that the modeler must justify by other means before drawing the conclusion.

\paragraph{Limitations.} The analysis is restricted to deterministic autonomous systems; stochastic, nonautonomous, delay, age-structured and spatial models are outside its scope. The mechanistic parameters $\theta$ are treated as known; joint estimation of $\theta$ and the closure raises questions of parametric identifiability and of entanglement between the two, visible already in \eqref{eq:R0seir}, where only the product $\beta^\ast\Lambda/\mu$ enters, that are not addressed here and for which the functional identifiability framework of \cite{LomanBaker2025,LomanBrowningBaker2026} is the natural tool. The identifiability results concern indistinguishability of exact trajectories on a set (Definition \ref{def:orbitid}); no statistical observation model, noise level or finite sample is involved, so they describe a structural limitation that persists with unlimited noiseless data on the visited region and say nothing about practical identifiability, which can only be worse. Global stability of the endemic equilibrium is established only for the separable subclass $\Ecal_{\mathrm{sep}}$, a genuine structural restriction; models with nonmonotone or multivariate closures may sustain oscillations or multiple endemic states, as the bifurcation analysis shows. The numerical section illustrates the constructions of the theory; it contains no trained model, and the frequency with which training produces closures of the kinds exhibited there is an empirical question left open.

\section{Conclusion}\label{sec:conclusion}

For universal differential equations whose closures parameterize intercompartmental fluxes and whose force of infection factorizes through the infected compartments, the qualitative theory of compartmental epidemic models transfers: positivity, dissipativity and invariance of the disease-free set hold for every weight vector; the reproduction number has a closed form; the disease-free and endemic equilibria, persistence and the bifurcation at threshold are analyzed with the classical tools under hypotheses that become architectural constraints or certifiable properties of the trained closure. The results that are specific to learned closures concern what the data determine. The reproduction number depends on the learned rates only through their values at the disease-free equilibrium, and, in the neural SEIR model, the bifurcation type additionally on the first derivatives of the transmission rate there; closures that coincide on the region visited by an observed trajectory are indistinguishable from data supported on that region, whatever is observed, and may differ in threshold and in bifurcation type. The uncertainty in the threshold is bounded by the Lipschitz constant of the closure times the distance of the data from the disease-free state, and monotone architectures give one-sided bounds. Agreement with trajectory data therefore does not by itself identify the invariants on which epidemiological interpretation rests; the architecture, the regularization and the data near the disease-free state do. Extending the framework to stochastic, nonautonomous, spatial and age-structured models, to joint inference of mechanistic parameters and closures, and to statistical rather than structural identifiability are the natural next steps.


\begin{thebibliography}{99}
\raggedright

\bibitem{Alsammani2026Tutorial} A. Alsammani, A tutorial on symbolic structural identifiability analysis of ODE models in Julia, \emph{Bulletin of Mathematical Biology} \textbf{88} (2026), 212. \doi{10.1007/s11538-026-01779-0}

\bibitem{AndersonMay1991} R.\,M. Anderson and R.\,M. May, \emph{Infectious Diseases of Humans: Dynamics and Control}, Oxford University Press, Oxford, 1991.

\bibitem{BermanPlemmons1994} A. Berman and R.\,J. Plemmons, \emph{Nonnegative Matrices in the Mathematical Sciences}, Classics in Applied Mathematics 9, SIAM, Philadelphia, 1994. \doi{10.1137/1.9781611971262}

\bibitem{BlowerDowlatabadi1994} S.\,M. Blower and H. Dowlatabadi, Sensitivity and uncertainty analysis of complex models of disease transmission: an HIV model, as an example, \emph{International Statistical Review} \textbf{62} (1994), 229--243. \doi{10.2307/1403510}

\bibitem{Brauer2004} F. Brauer, Backward bifurcations in simple vaccination models, \emph{Journal of Mathematical Analysis and Applications} \textbf{298} (2004), 418--431. \doi{10.1016/j.jmaa.2004.05.045}

\bibitem{BrauerCastilloChavezFeng2019} F. Brauer, C. Castillo-Chavez and Z. Feng, \emph{Mathematical Models in Epidemiology}, Texts in Applied Mathematics 69, Springer, New York, 2019. \doi{10.1007/978-1-4939-9828-9}

\bibitem{BruntonProctorKutz2016} S.\,L. Brunton, J.\,L. Proctor and J.\,N. Kutz, Discovering governing equations from data by sparse identification of nonlinear dynamical systems, \emph{Proceedings of the National Academy of Sciences} \textbf{113} (2016), 3932--3937. \doi{10.1073/pnas.1517384113}

\bibitem{CapassoSerio1978} V. Capasso and G. Serio, A generalization of the Kermack-McKendrick deterministic epidemic model, \emph{Mathematical Biosciences} \textbf{42} (1978), 43--61. \doi{10.1016/0025-5564(78)90006-8}

\bibitem{Carr1981} J. Carr, \emph{Applications of Centre Manifold Theory}, Applied Mathematical Sciences 35, Springer, New York, 1981. \doi{10.1007/978-1-4612-5929-9}

\bibitem{CastilloChavezSong2004} C. Castillo-Chavez and B. Song, Dynamical models of tuberculosis and their applications, \emph{Mathematical Biosciences and Engineering} \textbf{1} (2004), 361--404. \doi{10.3934/mbe.2004.1.361}

\bibitem{Chen2018NeuralODE} R.\,T.\,Q. Chen, Y. Rubanova, J. Bettencourt and D. Duvenaud, Neural ordinary differential equations, in \emph{Advances in Neural Information Processing Systems 31}, 2018, 6571--6583. arXiv:1806.07366.

\bibitem{Chitnis2008} N. Chitnis, J.\,M. Hyman and J.\,M. Cushing, Determining important parameters in the spread of malaria through the sensitivity analysis of a mathematical model, \emph{Bulletin of Mathematical Biology} \textbf{70} (2008), 1272--1296. \doi{10.1007/s11538-008-9299-0}

\bibitem{Cybenko1989} G. Cybenko, Approximation by superpositions of a sigmoidal function, \emph{Mathematics of Control, Signals and Systems} \textbf{2} (1989), 303--314. \doi{10.1007/BF02551274}

\bibitem{Dandekar2020} R. Dandekar, C. Rackauckas and G. Barbastathis, A machine learning aided global diagnostic and comparative tool to assess effect of quarantine control in COVID-19 spread, \emph{Patterns} \textbf{1} (2020), 100145. \doi{10.1016/j.patter.2020.100145}

\bibitem{DeRyckMishra2022} T. De Ryck and S. Mishra, Error analysis for physics informed neural networks (PINNs) approximating Kolmogorov PDEs, \emph{Advances in Computational Mathematics} \textbf{48} (2022), 79. \doi{10.1007/s10444-022-09985-9}

\bibitem{DiekmannHeesterbeekMetz1990} O. Diekmann, J.\,A.\,P. Heesterbeek and J.\,A.\,J. Metz, On the definition and the computation of the basic reproduction ratio $R_0$ in models for infectious diseases in heterogeneous populations, \emph{Journal of Mathematical Biology} \textbf{28} (1990), 365--382. \doi{10.1007/BF00178324}

\bibitem{DiekmannHeesterbeekRoberts2010} O. Diekmann, J.\,A.\,P. Heesterbeek and M.\,G. Roberts, The construction of next generation matrices for compartmental epidemic models, \emph{Journal of the Royal Society Interface} \textbf{7} (2010), 873--885. \doi{10.1098/rsif.2009.0386}

\bibitem{Dupont2019} E. Dupont, A. Doucet and Y.\,W. Teh, Augmented neural ODEs, in \emph{Advances in Neural Information Processing Systems 32}, 2019, 3140--3150. arXiv:1904.01681.

\bibitem{Dushoff1998} J. Dushoff, W. Huang and C. Castillo-Chavez, Backwards bifurcations and catastrophe in simple models of fatal diseases, \emph{Journal of Mathematical Biology} \textbf{36} (1998), 227--248. \doi{10.1007/s002850050099}

\bibitem{WeinanE2017} W. E, A proposal on machine learning via dynamical systems, \emph{Communications in Mathematics and Statistics} \textbf{5} (2017), 1--11. \doi{10.1007/s40304-017-0103-z}

\bibitem{FreedmanRuanTang1994} H.\,I. Freedman, S. Ruan and M. Tang, Uniform persistence and flows near a closed positively invariant set, \emph{Journal of Dynamics and Differential Equations} \textbf{6} (1994), 583--600. \doi{10.1007/BF02218848}

\bibitem{MohammedAlsammani2026}
M. Mohammed and A. Alsammani,
Long-term coexistence of epidemics and risk awareness:
impacts of adaptive human response and fatigue,
\emph{BIOMATH} (2026).
\doi{10.55630/j.biomath.2026.10.234}

\bibitem{GuckenheimerHolmes1983} J. Guckenheimer and P. Holmes, \emph{Nonlinear Oscillations, Dynamical Systems, and Bifurcations of Vector Fields}, Applied Mathematical Sciences 42, Springer, New York, 1983. \doi{10.1007/978-1-4612-1140-2}

\bibitem{HaberRuthotto2017} E. Haber and L. Ruthotto, Stable architectures for deep neural networks, \emph{Inverse Problems} \textbf{34} (2018), 014004. \doi{10.1088/1361-6420/aa9a90}

\bibitem{HadelerVandenDriessche1997} K.\,P. Hadeler and P. van den Driessche, Backward bifurcation in epidemic control, \emph{Mathematical Biosciences} \textbf{146} (1997), 15--35. \doi{10.1016/S0025-5564(97)00027-8}

\bibitem{HaleWaltman1989} J.\,K. Hale and P. Waltman, Persistence in infinite dimensional systems, \emph{SIAM Journal on Mathematical Analysis} \textbf{20} (1989), 388--395. \doi{10.1137/0520025}

\bibitem{Hartman2002} P. Hartman, \emph{Ordinary Differential Equations}, 2nd edition, Classics in Applied Mathematics 38, SIAM, Philadelphia, 2002. \doi{10.1137/1.9780898719222}

\bibitem{Hethcote2000} H.\,W. Hethcote, The mathematics of infectious diseases, \emph{SIAM Review} \textbf{42} (2000), 599--653. \doi{10.1137/S0036144500371907}

\bibitem{HornJohnson2013} R.\,A. Horn and C.\,R. Johnson, \emph{Matrix Analysis}, 2nd edition, Cambridge University Press, Cambridge, 2013. \doi{10.1017/CBO9781139020411}

\bibitem{Hornik1991} K. Hornik, Approximation capabilities of multilayer feedforward networks, \emph{Neural Networks} \textbf{4} (1991), 251--257. \doi{10.1016/0893-6080(91)90009-T}

\bibitem{HornikStinchcombeWhite1990} K. Hornik, M. Stinchcombe and H. White, Universal approximation of an unknown mapping and its derivatives using multilayer feedforward networks, \emph{Neural Networks} \textbf{3} (1990), 551--560. \doi{10.1016/0893-6080(90)90005-6}

\bibitem{Karniadakis2021} G.\,E. Karniadakis, I.\,G. Kevrekidis, L. Lu, P. Perdikaris, S. Wang and L. Yang, Physics informed machine learning, \emph{Nature Reviews Physics} \textbf{3} (2021), 422--440. \doi{10.1038/s42254-021-00314-5}

\bibitem{KermackMcKendrick1927} W.\,O. Kermack and A.\,G. McKendrick, A contribution to the mathematical theory of epidemics, \emph{Proceedings of the Royal Society of London A} \textbf{115} (1927), 700--721. \doi{10.1098/rspa.1927.0118}

\bibitem{Kidger2022} P. Kidger, On neural differential equations, PhD thesis, University of Oxford, 2022. arXiv:2202.02435.

\bibitem{Korobeinikov2007} A. Korobeinikov, Global properties of infectious disease models with nonlinear incidence, \emph{Bulletin of Mathematical Biology} \textbf{69} (2007), 1871--1886. \doi{10.1007/s11538-007-9196-y}

\bibitem{KorobeinikovMaini2004} A. Korobeinikov and P.\,K. Maini, A Lyapunov function and global properties for SIR and SEIR epidemiological models with nonlinear incidence, \emph{Mathematical Biosciences and Engineering} \textbf{1} (2004), 57--60. \doi{10.3934/mbe.2004.1.57}

\bibitem{Lagaris1998} I.\,E. Lagaris, A. Likas and D.\,I. Fotiadis, Artificial neural networks for solving ordinary and partial differential equations, \emph{IEEE Transactions on Neural Networks} \textbf{9} (1998), 987--1000. \doi{10.1109/72.712178}

\bibitem{LaSalle1976} J.\,P. LaSalle, \emph{The Stability of Dynamical Systems}, CBMS-NSF Regional Conference Series in Applied Mathematics 25, SIAM, Philadelphia, 1976. \doi{10.1137/1.9781611970432}

\bibitem{LiMuldowney1996} M.\,Y. Li and J.\,S. Muldowney, A geometric approach to global stability problems, \emph{SIAM Journal on Mathematical Analysis} \textbf{27} (1996), 1070--1083. \doi{10.1137/S0036141094266449}

\bibitem{Linka2022} K. Linka, A. Sch\"afer, X. Meng, Z. Zou, G.\,E. Karniadakis and E. Kuhl, Bayesian physics informed neural networks for real world nonlinear dynamical systems, \emph{Computer Methods in Applied Mechanics and Engineering} \textbf{402} (2022), 115346. \doi{10.1016/j.cma.2022.115346}

\bibitem{LiuLevinIwasa1986} W.\,M. Liu, S.\,A. Levin and Y. Iwasa, Influence of nonlinear incidence rates upon the behavior of SIRS epidemiological models, \emph{Journal of Mathematical Biology} \textbf{23} (1986), 187--204. \doi{10.1007/BF00276956}

\bibitem{LomanBaker2025} T.\,E. Loman and R.\,E. Baker, Functional and parametric identifiability for universal differential equations applied to chemical reaction networks, preprint, 2025. arXiv:2510.14140.

\bibitem{LomanBrowningBaker2026} T.\,E. Loman, A.\,P. Browning and R.\,E. Baker, Structural functional identifiability and model discovery in differential equation models, preprint, 2026. arXiv:2606.30289.

\bibitem{Marino2008} S. Marino, I.\,B. Hogue, C.\,J. Ray and D.\,E. Kirschner, A methodology for performing global uncertainty and sensitivity analysis in systems biology, \emph{Journal of Theoretical Biology} \textbf{254} (2008), 178--196. \doi{10.1016/j.jtbi.2008.04.011}

\bibitem{Martcheva2015} M. Martcheva, \emph{An Introduction to Mathematical Epidemiology}, Texts in Applied Mathematics 61, Springer, New York, 2015. \doi{10.1007/978-1-4899-7612-3}

\bibitem{Massaroli2020} S. Massaroli, M. Poli, J. Park, A. Yamashita and H. Asama, Dissecting neural ODEs, in \emph{Advances in Neural Information Processing Systems 33}, 2020, 3952--3963. arXiv:2002.08071.

\bibitem{McKay1979} M.\,D. McKay, R.\,J. Beckman and W.\,J. Conover, A comparison of three methods for selecting values of input variables in the analysis of output from a computer code, \emph{Technometrics} \textbf{21} (1979), 239--245. \doi{10.2307/1268522}

\bibitem{MishraMolinaro2022} S. Mishra and R. Molinaro, Estimates on the generalization error of physics informed neural networks for approximating PDEs, \emph{IMA Journal of Numerical Analysis} \textbf{43} (2023), 1--43. \doi{10.1093/imanum/drab093}

\bibitem{Pinkus1999} A. Pinkus, Approximation theory of the MLP model in neural networks, \emph{Acta Numerica} \textbf{8} (1999), 143--195. \doi{10.1017/S0962492900002919}

\bibitem{Rackauckas2020UDE} C. Rackauckas, Y. Ma, J. Martensen, C. Warner, K. Zubov, R. Supekar, D. Skinner, A. Ramadhan and A. Edelman, Universal differential equations for scientific machine learning, preprint, 2020. arXiv:2001.04385.

\bibitem{RackauckasNie2017} C. Rackauckas and Q. Nie, DifferentialEquations.jl: a performant and feature rich ecosystem for solving differential equations in Julia, \emph{Journal of Open Research Software} \textbf{5} (2017), 15. \doi{10.5334/jors.151}

\bibitem{RaissiPerdikarisKarniadakis2019}
M.~Raissi, P.~Perdikaris and G.~E.~Karniadakis, \newblock Physics-informed neural networks: A deep learning framework for solving forward and inverse problems involving nonlinear partial differential equations, \newblock \emph{Journal of Computational Physics} \textbf{378} (2019), 686--707. \doi{10.1016/j.jcp.2018.10.045}

\bibitem{RuanWang2003} S. Ruan and W. Wang, Dynamical behavior of an epidemic model with a nonlinear incidence rate, \emph{Journal of Differential Equations} \textbf{188} (2003), 135--163. \doi{10.1016/S0022-0396(02)00089-X}

\bibitem{Shaier2022} S. Shaier, M. Raissi and P. Seshaiyer, Data driven approaches for predicting spread of infectious diseases through DINNs: disease informed neural networks, \emph{Letters in Biomathematics} \textbf{9} (2022), 71--105.

\bibitem{ShuaiVdD2013} Z. Shuai and P. van den Driessche, Global stability of infectious disease models using Lyapunov functions, \emph{SIAM Journal on Applied Mathematics} \textbf{73} (2013), 1513--1532. \doi{10.1137/120876642}

\bibitem{Smith1995} H.\,L. Smith, \emph{Monotone Dynamical Systems: An Introduction to the Theory of Competitive and Cooperative Systems}, Mathematical Surveys and Monographs 41, American Mathematical Society, Providence, 1995.

\bibitem{SmithThieme2011} H.\,L. Smith and H.\,R. Thieme, \emph{Dynamical Systems and Population Persistence}, Graduate Studies in Mathematics 118, American Mathematical Society, Providence, 2011.

\bibitem{Thieme1993} H.\,R. Thieme, Persistence under relaxed point dissipativity (with application to an endemic model), \emph{SIAM Journal on Mathematical Analysis} \textbf{24} (1993), 407--435. \doi{10.1137/0524026}

\bibitem{Thieme2003} H.\,R. Thieme, \emph{Mathematics in Population Biology}, Princeton University Press, Princeton, 2003.

\bibitem{vandenDriesscheWatmough2002} P. van den Driessche and J. Watmough, Reproduction numbers and sub-threshold endemic equilibria for compartmental models of disease transmission, \emph{Mathematical Biosciences} \textbf{180} (2002), 29--48. \doi{10.1016/S0025-5564(02)00108-6}

\bibitem{Virmaux2018} A. Virmaux and K. Scaman, Lipschitz regularity of deep neural networks: analysis and efficient estimation, in \emph{Advances in Neural Information Processing Systems 31}, 2018, 3835--3844. arXiv:1805.10965.

\bibitem{Walter1998} W. Walter, \emph{Ordinary Differential Equations}, Graduate Texts in Mathematics 182, Springer, New York, 1998. \doi{10.1007/978-1-4612-0601-9}

\bibitem{XiaoRuan2007} D. Xiao and S. Ruan, Global analysis of an epidemic model with nonmonotone incidence rate, \emph{Mathematical Biosciences} \textbf{208} (2007), 419--429. \doi{10.1016/j.mbs.2006.09.025}

\bibitem{Zhao2017} X.-Q. Zhao, \emph{Dynamical Systems in Population Biology}, 2nd edition, CMS Books in Mathematics, Springer, Cham, 2017. \doi{10.1007/978-3-319-56433-3}

\bibitem{Zinihi2026} A. Zinihi, M. Ehrhardt and M.\,R. Sidi Ammi, Structure-preserving neural ODEs via nonstandard finite difference discretization, preprint, 2026. arXiv:2607.10858.

\end{thebibliography}
\end{document}